\documentclass[12pt,twoside]{article}
\newcommand{\frakg}{\ensuremath{\mathfrak g}}
\newcommand{\frakig}{\ensuremath{\mathfrak g}}
\newcommand{\fraksl}{\ensuremath{\mathfrak s}{\mathfrak l}}
\newcommand{\frakso}{\ensuremath{\mathfrak s}{\mathfrak o}}
\newcommand{\fraksp}{\ensuremath{\mathfrak s}{\mathfrak p}}
\usepackage{amsfonts,a4,epsf}
\newtheorem{prop}{Proposition}[section]
\newtheorem{lem}[prop]{Lemma}
\newtheorem{thm}[prop]{Theorem}
\newtheorem{rmk}[prop]{Remark}
\newtheorem{ex}[prop]{Example}
\newtheorem{df}[prop]{Definition}
\newtheorem{dt}[prop]{Definition and Theorem}
\newtheorem{cor}[prop]{Corollary}
\newtheorem{fact}[prop]{Fact}
\newcommand{\g}{${\mathfrak g}_2$ }
\newcommand{\gf}{${\mathfrak g}_2$}
\newcommand{\D}{$\hat{\cal D}_3$ }

\newcommand{\C}{{\bf C}}
\newcommand{\R}{{\bf{R}}}
\newcommand{\topf}{${\cal T}_{tpf}$}
\newcommand{\htopf}{$\check{\cal T}_{tpf}$}
\newcommand{\aw}{word}
\newcommand{\aws}{words}
\newcommand{\I}{I_{(\mathfrak g_2,V)}}
\newcommand{\It}{ $\I$ }
\newcommand{\openbox}{\leavevmode
  \hbox to.77778em{
  \hfil\vrule
  \vbox to.675em{\hrule width.6em\vfil\hrule}
  \vrule\hfil}}
\newcommand{\schluss}{\vspace{-4ex}  \hfill \ensuremath{\openbox} 
\vspace{4ex}}
\renewcommand{\Box}{\schluss}
\newcommand{\pstext}[1]{\hspace{1mm}\raisebox{-0.5ex}{\epsfysize=2.5ex
\epsffile{#1.eps}}\hspace{1mm}}
\newcommand{\psdiag}[3]{\hspace{1mm}\raisebox{-#1mm}{\epsfysize#2mm
\epsffile{#3.eps}}\hspace{1mm}}
\newcommand{\psbild}[2]{\vspace{3mm}\\ \mbox{\epsfysize=#1mm\epsffile{#2.eps}}
\vspace{3mm}}
\newcommand{\pst}[1]{\psdiag{0}{3}{#1}}

\title{The skein relation for the $(\frakg_2,V)$-link invariant}
\author{Anna-Barbara Berger\footnote{The authors are partially supported by the
Schweizerische Nationalfonds.\protect\newline
 Mathematisches Institut, Sidlerstr. 5, 3012 Bern, Switzerland,\protect\newline {\tt bergerab@math-stat.unibe.ch,
stassen@math-stat.unibe.ch}}\hspace{0.2em} and Ines Stassen$^*$}\date{April 1998}

\begin{document}
\maketitle

\begin{abstract}
Pulling back the weight system associated with the exceptional Lie
algebra $\frakg_2$ by a modification of the universal Vassiliev-Kontsevich
invariant yields a link invariant; extending it to 3-nets, we derive a
recursive algorithm for its evaluation.
\end{abstract}

\setcounter{section}{-1}
\section{Introduction}
There is a well-known technique for the construction of Vassiliev link
invariants: define a weight system (i.e.\ a linear form on the space of chord
diagrams respecting certain relations)
on the basis of some Lie algebraic data and pull it back by
the universal Vassiliev-Kontsevich invariant. But unfortunately, the latter
is not known explicitly enough to allow direct evaluation of these link
invariants.

Efforts have been made to handle the universal Vassiliev-Kontsevich
invariant by considering only
``elementary'' parts of links into which any link may be cut.
This approach has been successful in so far as one
may hope to find skein relations for the link invariants coming from Lie
algebras---a skein relation being an equation implying a recursive algorithm
for the computation of a link invariant, for example the one that
determines the famous Jones polynomial:
$$t^2 P(\psdiag{2}{6}{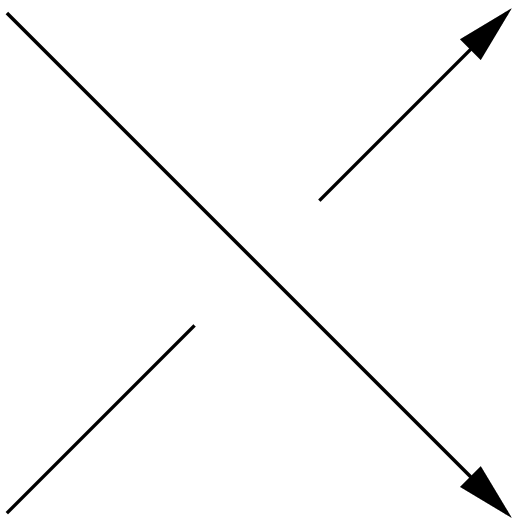})-t^{-2}P(\psdiag{2}{6}{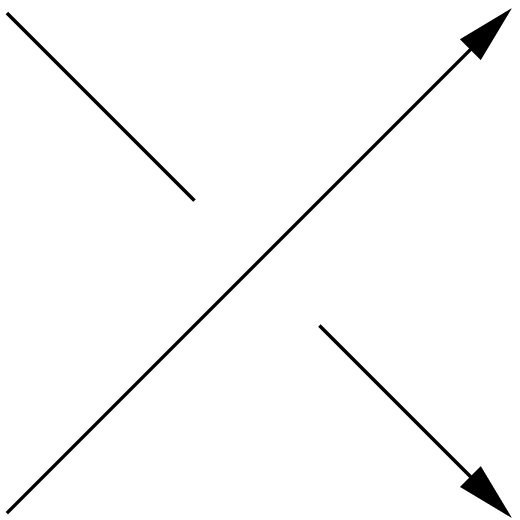})
= (t^{-1}-t)P(\psdiag{2}{6}{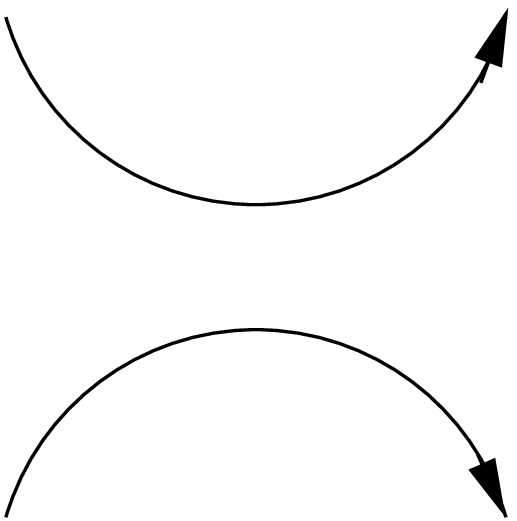}).$$
\medskip

It has been shown that the link invariants obtained from
the classical simple Lie algebras ${\fraksl}_n$, ${\frakso}_n$, and ${\fraksp}_n$
satisfy certain
versions of the skein relation of the HOMFLY polynomial (${\fraksl}_n$; see [LM 1])
resp.\ the Kauffman polynomial (${\frakso}_n$, ${\fraksp}_n$; see [LM 2]).
But what about the exceptional simple Lie algebras?
\medskip

In this paper, we deal with the case of the exceptional Lie algebra
$\frakg_2$. By means of a generalization of the notion of links---since we have
to introduce branchings, we call them 3-nets---we manage to establish a
skein relation for the $(\frakg_2,V)$-invariant, $V$ being the 7-dimensional
``standard'' representation of $\frakg_2$. As a by-product, we obtain an extension
of the ($\frakg_2,V$)-invariant to closed 3-nets. Kuperberg's skein relation
for the quantum $\frakg_2$-invariant (see [K]) turns out to be a special case of
ours; it is not too surprising that there is a connection between these
skein relations since the restrictions to knots of the two invariants
coincide according to a result of Piunikhin's in [P].

We expect that our method can be adapted to the case of the other
exceptional Lie algebras.
\medskip

To the reader not familiar with Lie theory, we recommend [H]
and [FH]. For an introduction to Vassiliev invariants and weight systems,
see [BN 1]; a more general definition of weight systems is given in [V],
section 6.
\medskip

{\bf Overview} over the categories and functors appearing in this paper:
\vspace{-2mm}\\
\begin{center}
$\I$ \\
\hspace*{9mm}\psdiag{0}{14}{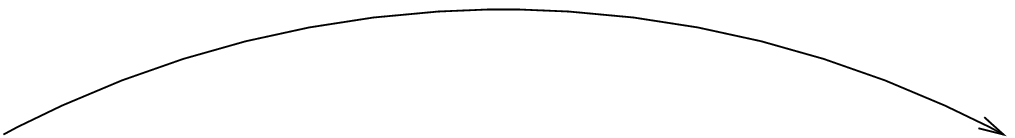} \vspace{-1.2cm}\\
\footnotesize the $(\frakg_2,V)$-invariant
\end{center}

\begin{tabular}{c c c c c}
\hspace{2.4cm} & \hspace{2.4cm} & \hspace{2.4cm} & \hspace{2.4cm}
& \hspace{2.4cm} \\
${\cal T}_{tpf}$ & $\stackrel{\hat Z_f}{\longrightarrow}$ & $\hat{\cal D}_3$ &
$\stackrel{\Psi}{\longrightarrow}$ & ${\cal C}(\frakg_2)$ \\
\end{tabular}

\begin{tabular}{p{2.4cm} p{2.4cm} p{2.4cm} p{2.4cm} p{2.4cm}}
\footnotesize the category of {\bf t}ri\-valent,
{\bf p}a\-ran\-the\-sized,
{\bf f}ramed tangles; eve\-ry 3-net corresponds to a morphism of this
category &
\footnotesize universal Vassi\-liev-Kontsevich in\-va\-riant &
\footnotesize the category of 3-diagrams (a generalization \newline of chord dia\-grams) &
\footnotesize $(\frakg_2,V)$-weight system &
\footnotesize a mo\-di\-fi\-ca\-tion of a sub\-ca\-te\-go\-ry of the ca\-te\-go\-ry of
re\-pre\-sen\-ta\-tions of $\frakg_2$
\end{tabular}
\vspace{3ex}

{\bf Acknowledgements} \hspace{0.3em} We would like to express our thanks to Ch.\ Riedtmann
for supervision, for encouragement, and for thoroughly reading this paper
and to P. Vogel for crucial suggestions and for patiently answering our
many questions; we also thank Ch.\ Blanchet for helpful discussions,
A. Beliakova for important references and K. Ott
for suggesting some linguistic improvement.




\section{3-nets and 3-tangles}

In this section, we will define 3-nets and 3-tangles. They are
generalizations of links and tangles. We will also construct the category
\topf \/ of trivalent  framed tangles, which will allow us to work in
a category-theoretical setting.\\

A 3-net will be something like a ``link with branchings''. To
describe the situation near a trivalent vertex (i.e.\ near a
branching point), we will need the following notion:\\ Let $B$ be
the unit ball in $\R^3$, i.e.\ $\{x\in
\R^3; |x|\leq 1\}$, together with the distinguished subset $T:= \{(t,0,0)|0\leq t\leq1\}\cup
\{(-\frac{1}{2},\frac{\sqrt{3}}{2},0)t|0\leq t\leq 1\}\cup
\{(-\frac{1}{2},-\frac{\sqrt{3}}{2},0)t|0\leq t\leq 1\}$.\\[1em]
\centerline{$B\quad =\quad$\psdiag{7}{21}{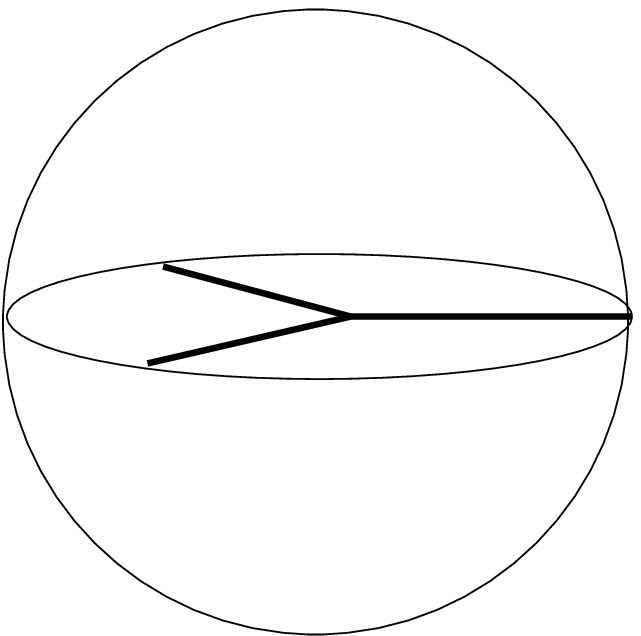}}

\begin{df}
A {\em framed 3-net} is a subset $N$ of \/ $\R^3$ with a finite
subset $\{t_1,\ldots ,t_n\}\subset N$ such that:
\begin{itemize}
\item[(i)] for every point $t_i$, there exists a
neighbourhood $U_i$ in $\R^3$ and a diffeomorphism
$f_i:U_i\rightarrow B$ with $f_i(N\cap U_i)= T$,
\item[(ii)] $\tilde{N}:= N \backslash \left(\bigcup_{i=1}^n f_i^{-1}(\{x\in B;
|x|<\frac{1}{2}\})\right)$ is an embedded  smooth  closed compact 1-dimensional
manifold,
\end{itemize}
together with:
\begin{itemize}
\item[(iii)] a  smooth vector field on $N$ that is nowhere tangent to $N$ (and in particullary nowhere zero)\footnote{Note that our framing is not a
framing in the classical sense.}.\\
\end{itemize}

  The points
$t_1,\ldots, t_n$ are called {\em trivalent vertices} of $N$;
boundary points $x$ of $\tilde{N}$ with $x\not\in U_i
\:(\forall i)$ are called {\em univalent vertices} of $N$.
\end{df}

Observe that 3-nets without trivalent and univalent vertices are
simply framed links. When we represent a framed 3-net $N$ by a diagram,
the framing of $N$ is given by the blackboard framing\footnote{i.e.\
the vector field assigned to $N$ consists of vectors pointing
upward }.
\newline The  3-net in the following figure is a 3-net
with 7 trivalent
and 3 univalent vertices.\\[1em]
\centerline{\psbild{25}{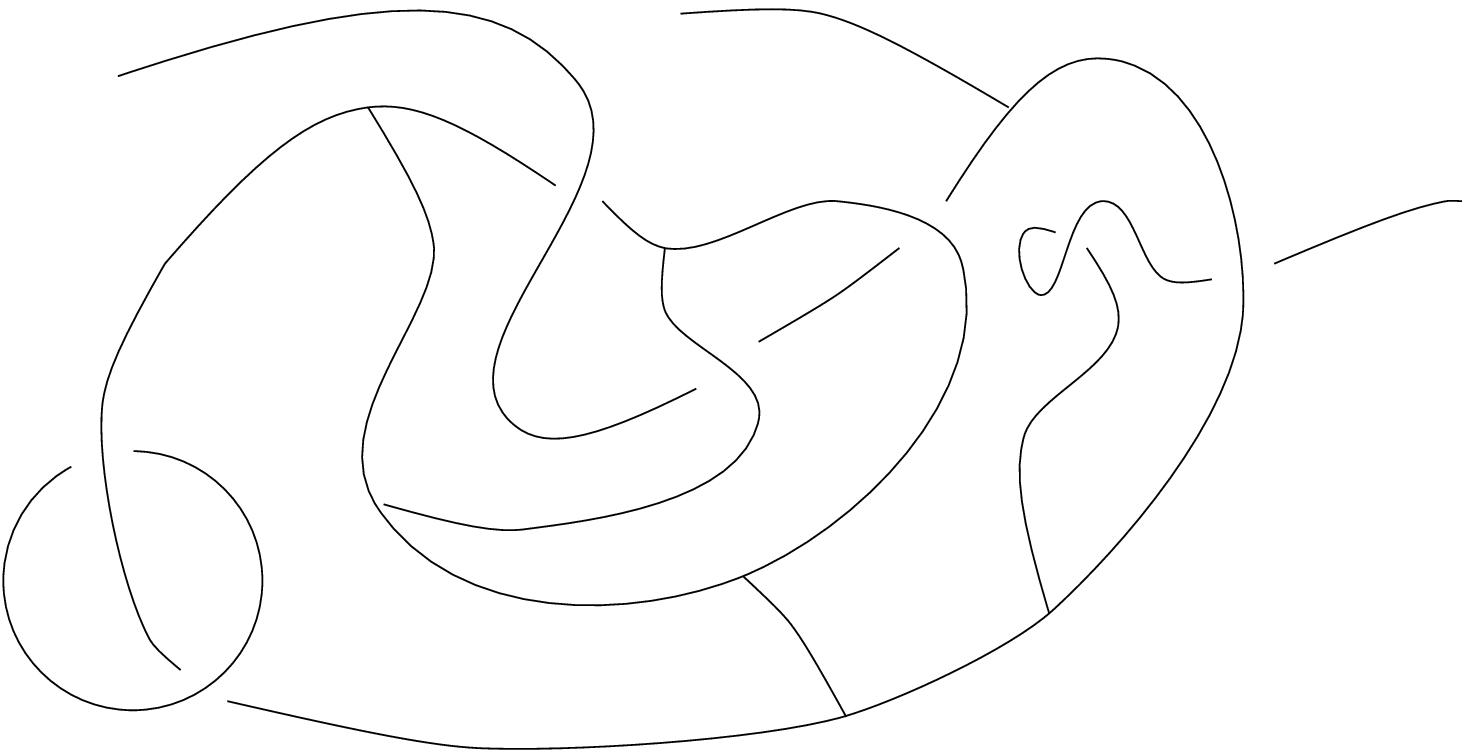}\quad\quad\quad\quad}

As all 3-nets in this paper will be framed 3-nets, we will usually
omit the word ``framed''.

\begin{df}\label{defeq}
Two 3-nets $N_1$ and $N_2$ are {\em equivalent} if $N_1$ can be
deformed into $N_2$ within the class of framed 3-nets by a smooth
1-parameter family of diffeomorphisms of $\R^3$.\newline A {\em
closed} 3-net is a 3-net without univalent vertices.\newline A
3-net is {\em planar} if it is equivalent to a 3-net $M$ with
$M\subset
\R^2\times \{0\}$ and the vector field assigned to $M$
consists of vectors of the form $(0,0,1)$.
\end{df}

\begin{rmk}
{\em The overview in the introduction intimates that we will obtain
an invariant of closed 3-nets that is composed of Vassiliev
invariants (see also section 6). Vassiliev invariants are usually
defined for oriented links, but we remind the reader that there is
a definition for unoriented links (see e.g.\ [St]):\newline
A link invariant\footnote{i.e.\ a function assigning to
each link an element of an abelian group (usually $\C$) that is
constant on the equivalence classes of framed links} $f$ is a {\it
Vassiliev invariant of type $m$} if for any link $L$, any diagram
$D(L)$ of $L$ and any subset $C$ of the set of crossings of $D(L)$
with cardinality greater than $m$ the following equation holds:
$$\sum_{X\subset C}(-1)^{|X|}f([D(L)_X])=0,$$
where $|X|$ is the cardinality of $X$, $D(L)_X$ is the link diagram
obtained form $D(L)$ by changing all the crossings in $X$, and
$[D(L)_X]$ is a link with diagram $D(L)_X$ (such that the framing
on the link is given by the blackboard framing of the diagram).
\newline
Of course, this definition can be extended to 3-nets. }\end{rmk}
\begin{df}
A {\em (framed) 3-tangle} is a framed 3-net $N$ with $N\subset
[0,1]\times
\R^2$ such that the points of $N$ lying in the planes $\{0\}\times
\R^2$ and $\{1\}\times \R^2$ are exactly the univalent vertices of
$N$ and these lie on one of the lines $\{0\}\times \R\times \{0\}$
and $\{1\}\times\R\times\{0\}$. Additionally,  we require that the
normal plane of $N$ in a univalent vertex $v$ is
parallel to the plane $\{0\}\times\R^2$, and the
vector field assigned to $N$ is $(0,0,1)$ in $v$.\end{df}

\begin{df}\label{deftangeq}
Two 3-tangles $T_1$ and $T_2$ are {\em equivalent} if one can be
deformed into the other within the class of 3-tangles by a smooth
1-parameter family of diffeomorphisms of $\R^3$.\newline A 3-tangle
is {\em planar} if it is equivalent to a 3-tangle $M$ with
$M\subset [0,1]\times\R\times \{0\}$ and the vector field assigned
to $M$ consists of vectors of the form $(0,0,1)$.
\end{df}

Now we will define the category of {\bf t}rivalent, {\bf
p}arenthesized,
{\bf f}ramed tangles \topf. It is an (unoriented) generalization of the category of
non-associative tangles in [BN 2]. \\

First, we will define non-associative words (which will be the objects of
\topf).
\begin{df}
A {\em non-assocative \aw} is a word $w$ in the alphabet
$\{\diamond, ), (\}$ such that $w$ is equal to the empty word,
$\diamond, (w_1)$ or $(w_1w_2)$ where $w_1, w_2$ are
non-associative
\aws.\newline
The {\em length} $l(w)$ of a non-associative word $w$ is the number of symbols
$\diamond$ in $w$.
\end{df}

\begin{ex}\quad\quad\quad\quad\quad $l\left(\diamond((\diamond(\diamond(\diamond\diamond)))\diamond)\right)
= 6$
\end{ex}

\begin{df}
Let \htopf \/ be the monoidal $\C$-category whose objects are
non-as\-so\-cia\-tive words (where the tensor product $w_1\otimes
w_2$ is $(w_1w_2)$ and the unit object is the empty word) and whose
morphisms are freely generated by the following morphisms:
\begin{itemize}
\item[(G1)] A morphism ${\mbox{\pstext{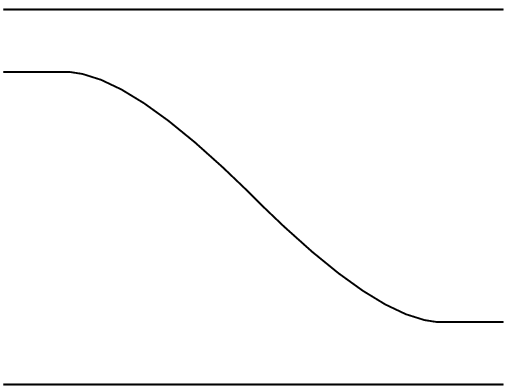}}}_{v,w,x}$ and a morphism
 ${\mbox{\pstext{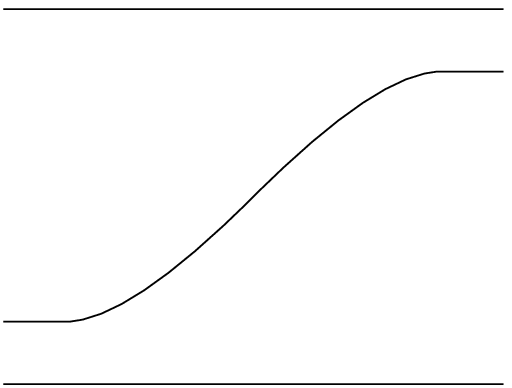}}}_{v,w,x}$ for each triple $(v,w,x)$ of
non-associative \aws: The sources of these morphisms are $((vw)x)$
and $(v(wx))$ and their targets are $(v(wx))$ and $((vw)x)$
respectively.\newline Graphically, we represent these morphisms as
the following examples indicate (the parenthesation of source and
target is encoded in the distances between the strands):
\end{itemize}
\begin{tabbing}
xxxxxxxxxxxxxxxxxxxxxxx \= xxxxxxxxxxxxxx \=  \kill
\>${\mbox{\pstext{phi}}}_{ (\diamond\diamond),\diamond,(\diamond\diamond)}\quad\quad$
\>$\leftrightarrow
\quad\quad\mbox{\psdiag{4}{12}{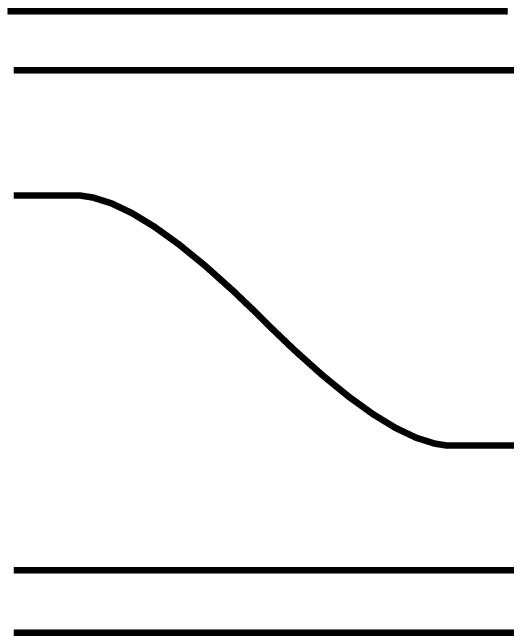}}$\\[2ex]
\>${\mbox{\pstext{phi-}}}_{\diamond,(\diamond\diamond),(\diamond\diamond)}\quad\quad$
\>$\leftrightarrow \quad\quad \mbox{\psdiag{4}{12}{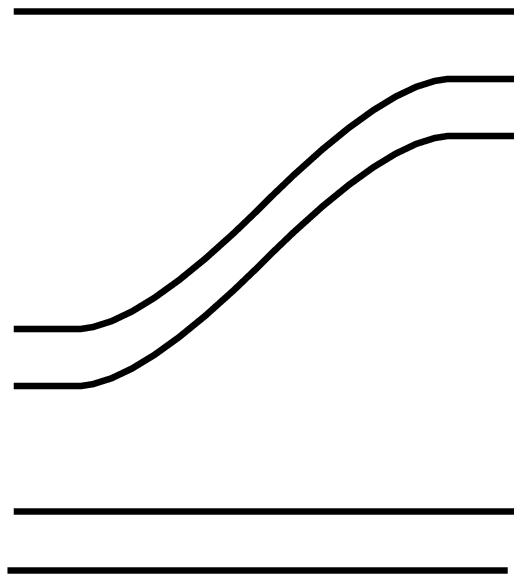}}$
\end{tabbing}
\begin{itemize}
\item[(G2)]
A morphism ${\mbox{\pstext{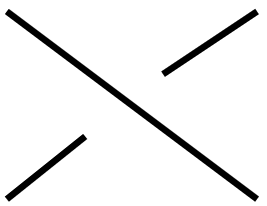}}}_{v,w}$ and a morphism
 ${\mbox{\pstext{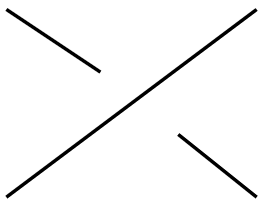}}_{v,w}}$ for each pair $(v,w)$ of
non-associative \aws: The source of these morphisms is
$(v,w)$  and their target is
$(w,v)$. We will depict them as follows:
\end{itemize}
\begin{tabbing}
xxxxxxxxxxxxxxxxxxxxxxx \= xxxxxxxxxxxxxx \=  \kill
\>${\mbox{\pstext{Zarg1}}}_{\diamond,((\diamond\diamond)\diamond)}\quad\quad$
\>$\leftrightarrow\quad\quad\mbox{\psdiag{3.5}{10.5}{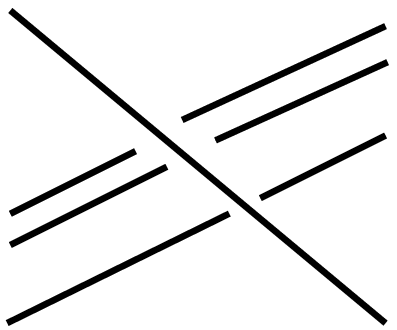}}$\\[2ex]
\>${\mbox{\pstext{Zarg5}}}_{(\diamond(\diamond\diamond)),(\diamond\diamond)}\quad\quad$
\>$\leftrightarrow\quad\quad\mbox{\psdiag{4}{12}{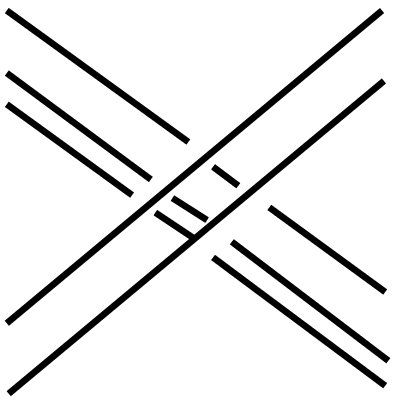}}$
\end{tabbing}
\begin{itemize}
\item[(G3)]A morphism ${\mbox{\pstext{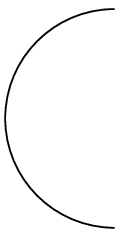}}}$ and a morphism \pstext{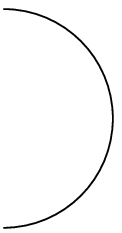}
with sources  $(\:)$ and $(\diamond\diamond)$ and targets $(\diamond\diamond)$ and
$(\:)$ respectively; graphically: \end{itemize}
\begin{tabbing}
xxxxxxxxxxxxxxxxxxxxxxx \= xxxxxxxxxxxxxx \=  \kill
\>${\mbox{\pstext{lbog}}}\quad\quad$\>$\leftrightarrow\quad\quad
\mbox{\psdiag{3}{9}{lbog}}$\\[2ex]
\>${\mbox{\pstext{rbog}}}\quad\quad$\>$\leftrightarrow\quad\quad
\mbox{\psdiag{3}{9}{rbog}}$
\end{tabbing}
\begin{itemize}
\item[(G4)]
A morphism ${\mbox{\pstext{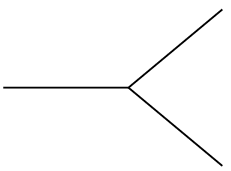}}}$ and a morphism ${\mbox{\pstext{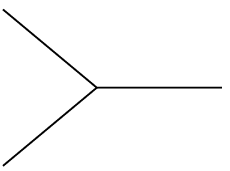}}}$ with
sources
$\diamond$ and $(\diamond\diamond)$ and targets $(\diamond\diamond)$ and $\diamond$
respectively; graphically:
\end{itemize}
\begin{tabbing}
xxxxxxxxxxxxxxxxxxxxxxx \= xxxxxxxxxxxxxx \=  \kill
\>${\mbox{\pstext{trir}}}\quad\quad$\>$\leftrightarrow\quad\quad
\mbox{\psdiag{3}{9}{trir}}$\\[2ex]
\>${\mbox{\pstext{tril}}}\quad\quad$\>$\leftrightarrow\quad\quad
\mbox{\psdiag{3}{9}{tril}}$
\end{tabbing}
The graphical representation of the tensor product of two morphisms
is obtained by putting the first above the second, the graphical
representation of the composition $S_2\circ S_1$ is obtained by
glueing the graphical representation of $S_2$ to the one of $S_1$
from the right:
\begin{tabbing}
xxxxxxxxxxxxxxxxxxxxxxx \= xxxxxxxxxxxxxx \=  \kill
\>$S_1 \otimes S_2\quad\quad $\>$\leftrightarrow \quad\quad
\mbox{\psdiag{5}{12}{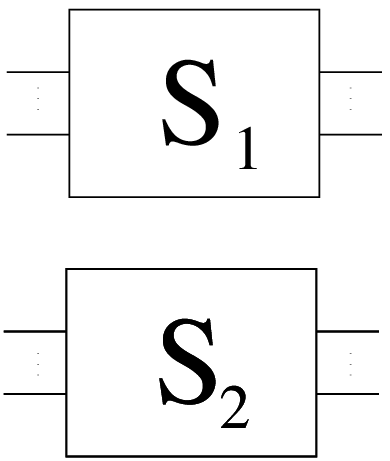}}$\\[3ex]
\>$S_2 \circ S_1\quad\quad $\>$\leftrightarrow \quad\quad
\mbox{\psdiag{2}{6}{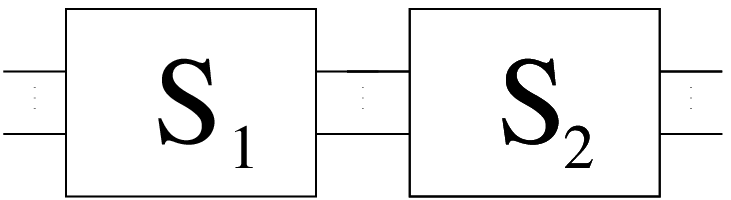}}$
\end{tabbing}
\end{df}

The graphical representations of the morphisms allow us to assign a
3-tangle to each morphism $M$ (namely a 3-tangle $T$ with the
graphical representation of $M$ as diagram, supplied with the
blackboard framing).

\begin{df}
The monoidal $\C$-category \topf \/ is the category whose
objects are the objects of \htopf \/ and whose morphisms are
the equivalence
classes of the morphisms of \htopf \/ under the following
equivalence relation: Two morphisms from $u$ to $w$ are equivalent
if they get assigned equivalent 3-tangles.
\end{df}

\begin{rmk}\label{lokrel} {\em The equivalence relation in the above definition can be
described locally. For morphisms generated by (G1)-(G3) this is done in [BN
2]. If we take the generators in  (G4) as well, we have to add the following
relations for any word $w$:

\begin{tabbing}
\underline{Relation}\hspace{8cm}\= \underline{Graphical
representation}\\[1em]
$(id_\diamond\otimes\mbox{\pstext{rbog}})\,\pstext{phi}_{\diamond,\diamond,\diamond}\,({\mbox{\pstext{trir}}}\otimes
id_\diamond)
= \mbox{\pstext{tril}}$
\>\psdiag{2}{8}{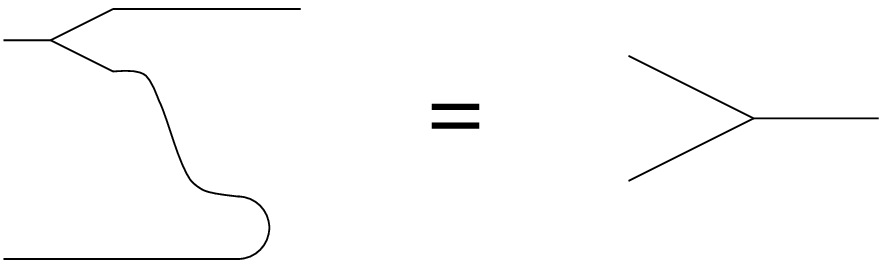}\\[1em]
$(id_\diamond\otimes
\mbox{\pstext{tril}})\,\pstext{phi}_{\diamond,\diamond,\diamond}\,({\mbox{\pstext{lbog}}}\otimes
id_\diamond)
={\mbox{\pstext{trir}}}$
\>\psdiag{2}{8}{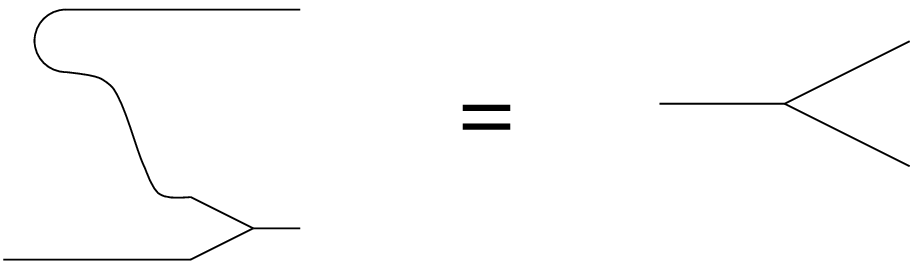}\\[1em]
$(id_w\otimes{\mbox{\pstext{tril}}}){\mbox{\pstext{Zarg1}}}_{(\diamond\diamond),w}
={\mbox{\pstext{Zarg1}}}_{\diamond,w}({\mbox{\pstext{tril}}}\otimes id_w)$
\>\psdiag{2}{6}{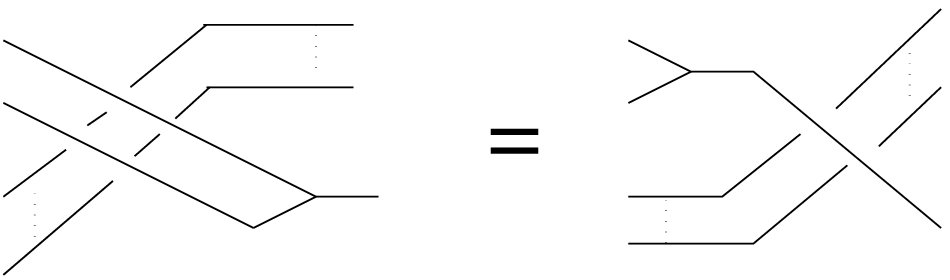}\\[1em]
$({\mbox{\pstext{tril}}}\otimes id_w){\mbox{\pstext{Zarg1}}}_{w,(\diamond\diamond)}
={\mbox{\pstext{Zarg1}}}_{w,\diamond}(id_w \otimes {\mbox{\pstext{tril}}})$
\>\psdiag{2}{6}{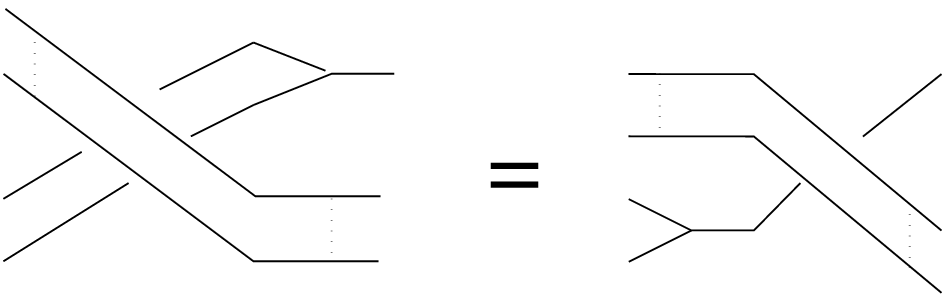}\\

\end{tabbing}}
\end{rmk}
\begin{rmk}
{\em Observe that any 3-tangle, and in particular any closed 3-net,
is equivalent to a 3-tangle assigned to a morphism of
\topf \/, and so the equivalence classes of 3-tangles correspond
exactly to a basis of the morphisms of \topf. One might
achieve this by taking much simpler generators (e.g.\
only generators without multiple strands),
 but then the local description of the equivalence relation given in remark
\ref{lokrel} would be more
complicated.}
\end{rmk}

\section{The universal Vassiliev-Kontsevich invariant extended to
3-tangles}

Now we want to extend and adapt the universal Vassiliev-Kontsevich invariant to the
3-tangles in ${\cal T}_{tpf}$. Since this invariant operates by taking a diagram
representing the given tangle as support and adding some chords,
we have to generalize the
notion of chord diagram and introduce trivalent vertices
in the support, too.

\begin{df}
A {\em 3-diagram} is a finite trivalent graph K (by which we understand
a graph with every
vertex being either univalent or trivalent or else bivalent and adjacent
to a loop) equipped with the following data:
\begin{itemize}
\item a colouring of the edges by \psdiag{0.3}{3}{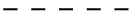} or
\psdiag{0.3}{3}{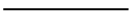} such that there is not a
vertex adjacent to two edges coloured by \psdiag{0.3}{3}{sstrich} and one coloured by
\psdiag{0.3}{3}{dstrich}
\item a colouring of the univalent vertices by $\,\circ\,$ or $\,\bullet\,$
according to whether the edge arriving there is coloured by \psdiag{0.3}{3}{sstrich} or
\psdiag{0.3}{3}{dstrich}
\item for every trivalent vertex x of K, a cyclic order of the edges
arriving at x.
\end{itemize}
The union of the edges coloured by
\psdiag{0.3}{3}{dstrich} is referred to as the {\em support} of the
diagram; the edges coloured by \psdiag{0.3}{3}{sstrich} are called {\em
chords}.\\
The {\em degree} of a 3-diagram is the number of trivalent
vertices adjacent to at least one chord \footnote{Note that for a 3-diagram without
univalent vertices adjacent to a chord, this is twice the classical
degree.}.
\end{df}

Usually, we describe the 3-diagrams by graphical representations
in the plane encoding
the information about the cyclic order near the trivalent vertices by
arranging the adjacent edges counterclockwise.

\begin{df}
The category ${\cal D}_3$ is a monoidal\/ {\bf C}-category whose morphisms are
linear combinations of
certain graphical representations of 3-diagrams. It is given by the
following data:
\\
{\bf objects:} {\em Obj}$({\cal D}_3) := \,\,^{\cdot}\hspace{-2.8mm}\bigcup_{n=0}^\infty\{\circ,\bullet\}^n$.
The tensor product on {\em Obj}$({\cal D}_3)$ is the juxtaposition.
\\
{\bf generators:} The morphism spaces are generated by:
\\  \psbild{36}{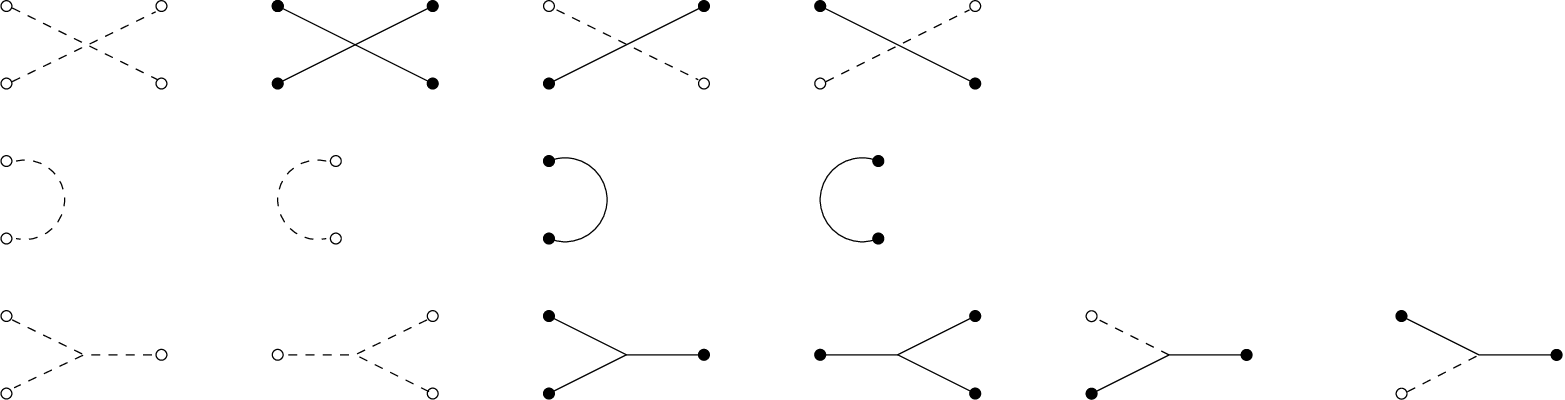}
\\ the source (resp.\ the target) being denoted on the left-(resp.\
right-)hand side from top to bottom
\footnote{The apparent 4-valent vertices are no vertices at all - they are just
crossings of two edges (there is no need to say that one of them passes
over the other).}.
\\ The tensor product of two morphisms is obtained by putting the first above the
second, the composition by glueing together the corresponding entries of
the target of the first and the source of the second.
\\
{\bf relations:} Of course, different graphical representations of
isomorphic 3-dia\-grams are to represent the same morphism; in addition, we
impose the following re\-la\-tions
\footnote{The reader familiar with Bar-Natan's way of defining weight systems
should pay attention to the sign in our relation {\em (AS3)}.}:
\begin{tabbing}
(AS1) xxxxxxx \= \kill
(AS1) \> $\psdiag{3}{9}{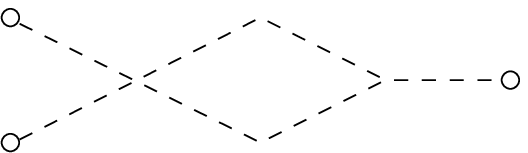} \qquad = \qquad - \psdiag{3}{9}{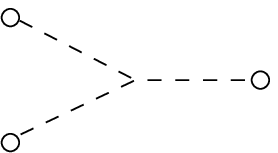}$ \\[1ex]
(AS2) \> $\psdiag{3}{9}{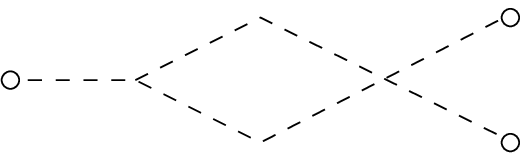} \qquad = \qquad - \psdiag{3}{9}{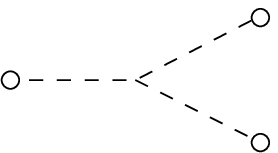}$ \\[1ex]
(AS3) \> $\psdiag{3}{9}{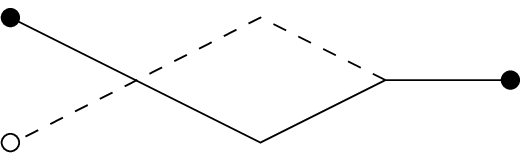} \qquad = \qquad - \psdiag{3}{9}{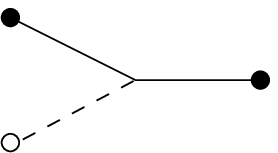} $\\[1ex]
(IHX1) \> $\psdiag{4}{12}{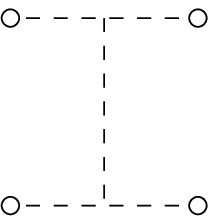} \qquad = \qquad \psdiag{4}{12}{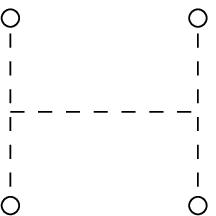}  -  \psdiag{4}{12}{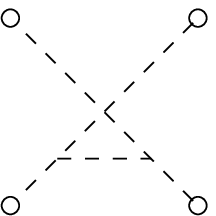}$
\\[1ex]
(IHX2) \>$\psdiag{4}{12}{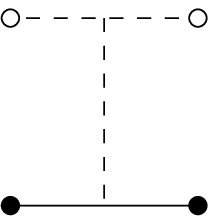} \qquad = \qquad \psdiag{4}{12}{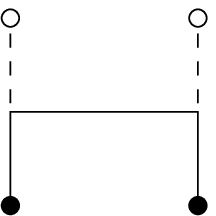}  -  \psdiag{4}{12}{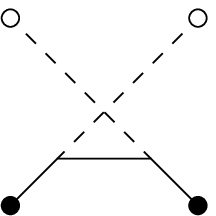}$
\\[1ex]
(IHX3) \> $\psdiag{4}{12}{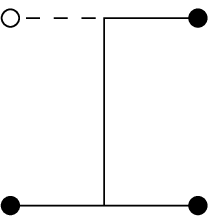}\qquad = \qquad \psdiag{4}{12}{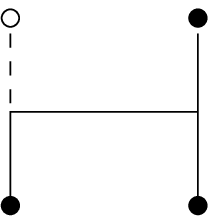}  -  \psdiag{4}{12}{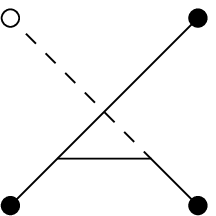}$
\end{tabbing}
\end{df}

Obviously, the morphisms involved in the relations described above can be
composed of the generators of ${\cal D}_3$; the relation {\em (IHX2)} for
example can be written as follows:
\begin{center}
\psdiag{6}{17}{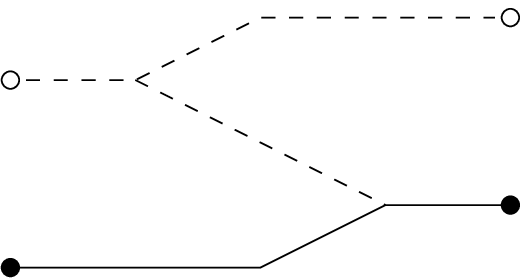} = \psdiag{7.5}{21.5}{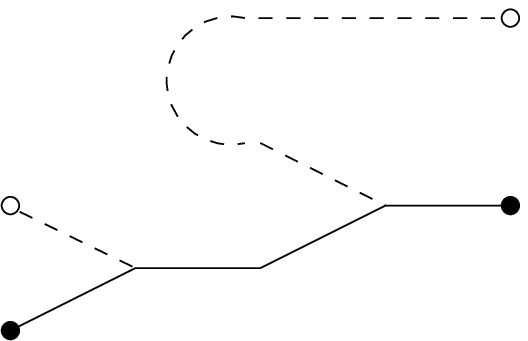} - \psdiag{9}{26}{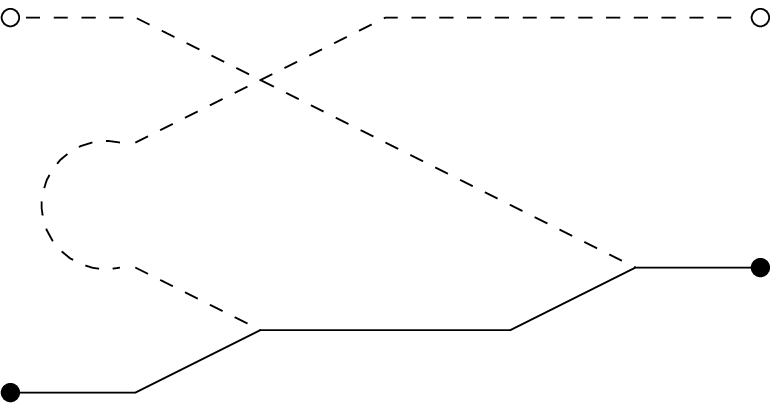}
\end{center}
But for obvious reasons, we refrain from doing so.

\begin{rmk}
{\em The relations {\em (AS1)}, {\em (AS2)}, and {\em (IHX1)} allow
introducing trivalent vertices adjacent to three chords in a consistent
way, which are convenient e.g.\ for the formulation of {\em (IHX2)}.
The relations {\em (AS3)}, {\em (IHX2)}, and {\em (IHX3)} are required for
the existence of the universal Vassiliev-Kontsevich invariant. See proof of
\ref{VK} for {\em (AS3)}
and proof of theorem 1 (1) in [BN 1] for {\em (IHX2)} (the 4T-relation);
{\em (IHX3)} reflects a similar situation near a trivalent vertex. }
\end{rmk}

\begin{df}
Let $\hat{\cal D}_3$ be the completion of the (graded) category $\hat{\cal
D}$.
\end{df}

For convenience of notation, we define a functor $\Delta$ from a
specialized version ${\cal D}^*_3$ of the category $\hat{\cal D}_3$ into
$\hat{\cal D}_3$.

\begin{df}
Let ${\cal D}^*_3$ be the category whose objects and morphisms are those of
$\hat{\cal D}_3$ together with some extra information:  in the morphisms, some
connected components of the support containing no trivalent vertices
adjacent to
three edges coloured by \psdiag{0.3}{3}{dstrich} may be labelled by replacing
the adjacent components of the source and/or the
target by $*$.
\\ The composition of morphisms is to respect the labelling.
\\ The relations of ${\cal D}^*_3$ are those of $\hat{\cal D}_3$ with any
possible labelling.
\end{df}

If $D$ is an arbitrary morphism in $\hat{\cal D}_3$, we denote by $D_k$ the
morphism of ${\cal D}^*_3$ that is $D$ with the component departing from
the $k$-th entry of the source labelled.

\begin{df}
Let $\Delta$ be the monoidal functor ${\cal D}^*_3 \rightarrow \hat{\cal D}_3$
doubling the labelled component and taking the sum over all possible ways
of lifting arriving chords to the new components of the support; thus,
$\Delta$ is given by: \vspace{3mm}
\begin{tabbing}
$\Delta$(\psdiag{2}{6}{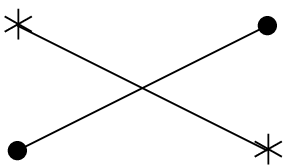})\hspace{1mm} \= := \hspace{1mm}
\=\psdiag{5}{12}{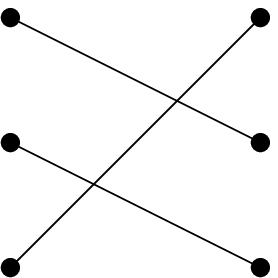} \hspace{3mm} \=
$\Delta$(\psdiag{2}{6}{deltarg1})\hspace{1mm} \= := \hspace{1mm}
\=\psdiag{5}{12}{deltim1} \hspace{3mm} \=
$\Delta$(\psdiag{2}{6}{deltarg1})\hspace{1mm} \= := \hspace{1mm}
\=\psdiag{5}{12}{deltim1} \kill
$\Delta(*)$ \> := \> $\bullet\enspace\bullet$
\> $\Delta(\bullet)$ \> := \> $\bullet$
\> $\Delta(\circ)$ \> := \> $\circ$ \\
[3ex]
$\Delta$(\psdiag{2}{6}{deltarg1}) \> := \> \psdiag{5}{12}{deltim1}
\> $\Delta$(\psdiag{2}{6}{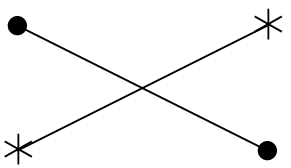}) \> := \> \psdiag{5}{12}{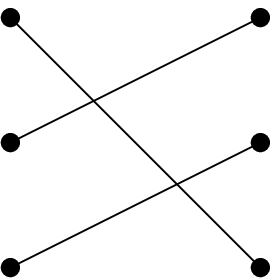}
\> $\Delta$(\psdiag{2}{6}{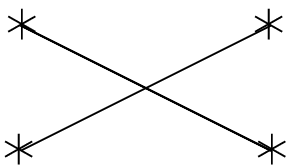}) \> := \> \psdiag{8}{18}{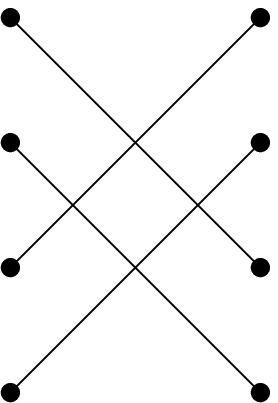} \\
[2ex]
$\Delta$(\psdiag{2}{6}{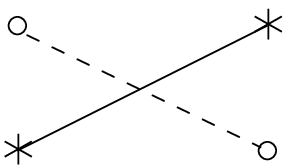}) \> := \> \psdiag{5}{12}{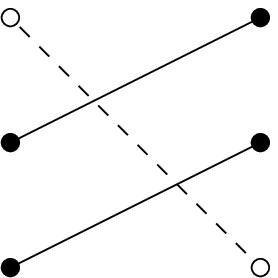}
\> $\Delta$(\psdiag{2}{6}{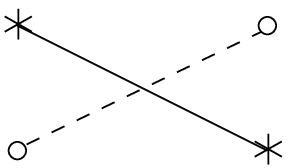}) \> := \> \psdiag{5}{12}{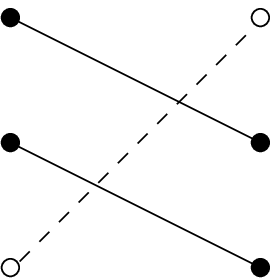} \\
[2ex]
$\Delta$(\psdiag{2}{6}{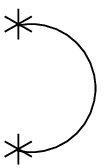}) \> := \> \psdiag{8}{18}{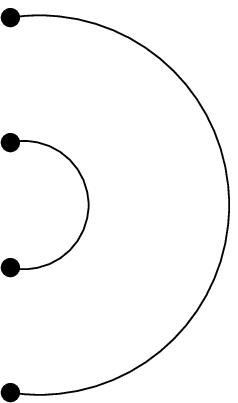}
\>$\Delta$(\psdiag{2}{6}{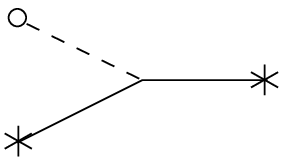}) \> := \> \psdiag{5}{12}{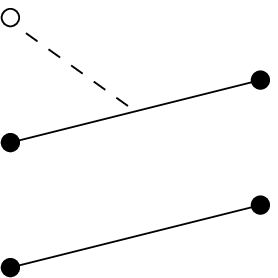}
+ \psdiag{5}{12}{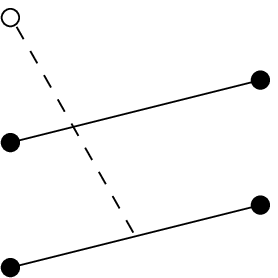} \\
[2ex]
$\Delta$(\psdiag{2}{6}{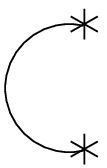}) \> := \> \psdiag{8}{18}{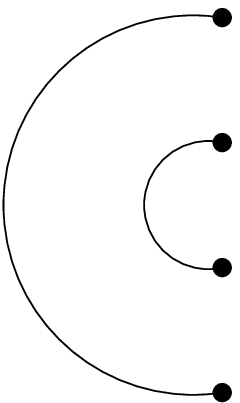}
\> $\Delta$(\psdiag{2}{6}{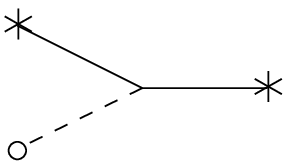}) \> := \> \psdiag{5}{12}{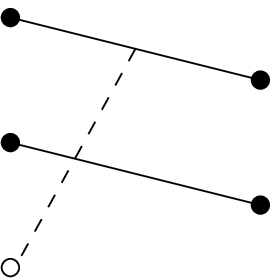}
+ \psdiag{5}{12}{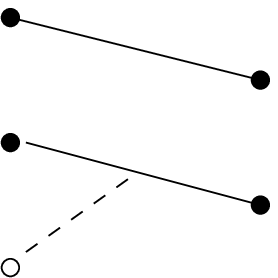}
\end{tabbing}

All the other generators (without any labelling) are not affected by
$\Delta$.
\end{df}

As an immediate consequence of the relation {\em (IHX3)}, we obtain:

\begin{lem}\label{lemma2}
\psdiag{5}{9}{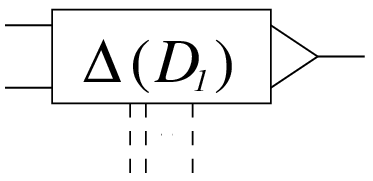} = \psdiag{5}{9}{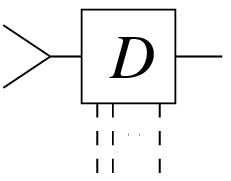},
\psdiag{5}{9}{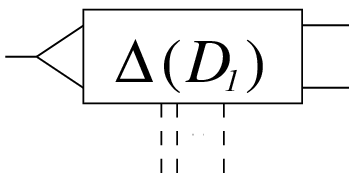} = \psdiag{5}{9}{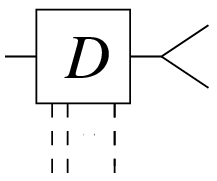}
\\ for any fitting 3-diagram D. \hfill $\Box$
\end{lem}

Finally, we get to the definition of the universal Vassiliev-Kontsevich
invariant:

\begin{dt}\label{VK}
The following assignments define a monoidal
functor $\hat Z_f : {\cal T}_{tpf} \to \hat{\cal D}_3$, the
{\em (unoriented) universal Vassiliev-Kontsevich in\-va\-riant}:
\vspace{2mm}\\
$\hat Z_f(u):=\bullet\,^n$, where $u$ is a non-associative
word of length $n$.
\begin{tabbing}
$\hat Z_f$(\psdiag{3}{9}{Zarg1})xx \= := xx\=\psdiag{3}{9}{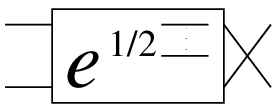} xxxxxx
\= $\hat Z_f$(\psdiag{3}{9}{Zarg1})xx \= := xx\=\psdiag{3}{9}{Zim1}
\kill
$\hat Z_f$(\pstext{Zarg1}) \> := \> \psdiag{2}{6}{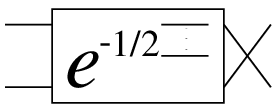}
\> $\hat Z_f$(\pstext{Zarg5}) \> := \> \psdiag{2}{6}{Zim1} \\
[2ex]
$\hat Z_f$(\pstext{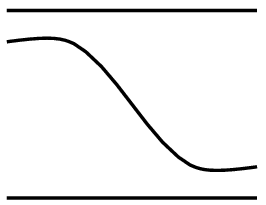}) \> := \> \psdiag{2}{6}{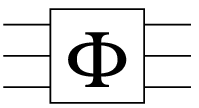}
\> $\hat Z_f$(\pstext{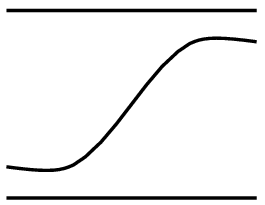}) \> := \> \psdiag{2}{6}{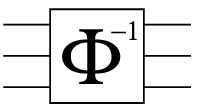} \\
[1ex]
$\hat Z_f$(\pstext{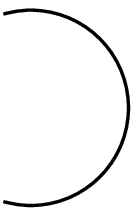}) \> := \> \psdiag{2}{6}{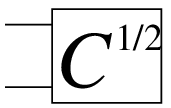}
\> $\hat Z_f$(\pstext{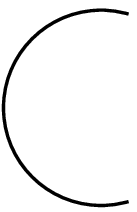}) \> := \> \psdiag{2}{6}{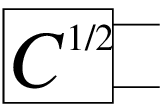} \\
[1ex]
$\hat Z_f$(\pstext{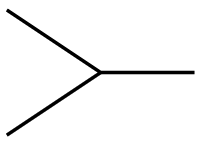}) \> := \> $\hat r\cdot$\psdiag{2}{6}{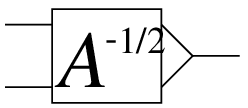}
\> $\hat Z_f$(\pstext{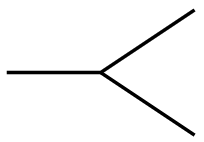}) \> := \> $\hat r\cdot$\psdiag{2}{6}{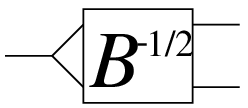} \\
\end{tabbing}
\begin{tabbing}
where xx \= \kill
where \> $e^{\pm \frac{1}{2}\psdiag{0}{2}{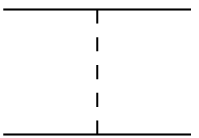}} := \sum_{n=0}^{\infty}$
$(\pm\frac{1}{2})^n\frac{1}{n!}\psdiag{1}{4}{hh}^n$ \\
      \> $\Phi$ is the Knizhnik-Zamolodchikov associator (for definition see [LM 3]) \\
      \> C := \pstext{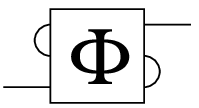} \\
      \> $\hat r \in \C \setminus \{ 0 \} $ \\
      \> A := \psdiag{3}{10}{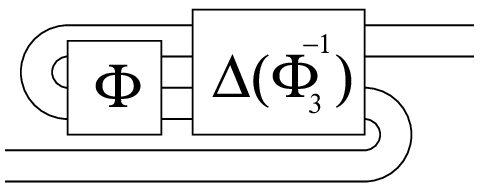}, B := \psdiag{3}{10}{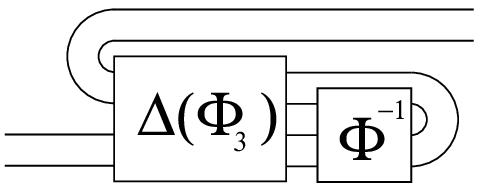}.
\end{tabbing}
For generators with multiple strands, i.e.\ $\pstext{Zarg1}_{u,v}$,
$\pstext{Zarg5}_{u,v}$, $\pstext{Zarg2}_{u,v,w}$ or $\pstext{Zarg6}_{u,v,w}$
for some non-associative words u, v, and w, the image
under $\hat Z_f$ is obtained by splitting the corresponding components of
$\hat Z_f(\pstext{Zarg1}_{\diamond,\diamond})$,
$\hat Z_f(\pstext{Zarg5}_{\diamond,\diamond})$,
$\hat Z_f(\pstext{Zarg2}_{\diamond,\diamond,\diamond})$, or
$\hat Z_f(\pstext{Zarg6}_{\diamond,\diamond,\diamond})$ respectively
by repeated
use of $\Delta$, e.g.\
$$\hat Z_f(\pstext{Zarg1}_{u,v}) =
\underbrace{\Delta\ldots ( \Delta}_{l(v)-1\, times}
(\underbrace{\Delta\ldots ( \Delta}_{l(u)-1\, times}
(\hat Z_f(\pstext{Zarg1})_1)_{\ldots 1})_{l(u)+1})_{\ldots l(u)+1}).$$
\end{dt}

\begin{rmk}
{\em As $\Phi$, A, B, C are formal power series in certain 3-diagrams
with degree 0-part 1, one may take
their inverses and square roots by substituting $x$ for their higher degree
parts and expand the corresponding function of $x$ in a Taylor series. }
\end{rmk}

\begin{rmk}
{\em Since the number of trivalent vertices in a 3-tangle is invariant, $\hat r$
may be chosen freely.}
\end{rmk}

{\bf Proof of \ref{VK}:} In section 1 of [MO], Murakami and Ohtsuki define the universal
Vassiliev-Kontsevich invariant for {\em oriented} 3-tangles and prove that
it is indeed an invariant. Their definition can be modified---without
destroying the invariance of the functor---as follows: Omit the signs
accounting for the orientation of the strands, and introduce the
antisymmetry relation {\em (AS3)} instead of the corresponding (implicit)
symmetry relation in the category of oriented 3-diagrams. In order to
obtain as degree 1-part of
$\hat Z_f(\pstext{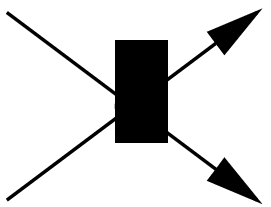}) := \hat Z_f(\pstext{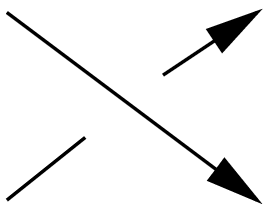}) - \hat Z_f(\pstext{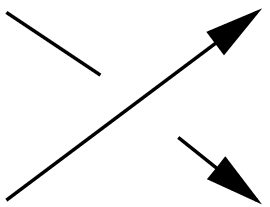})$
the diagram in which the double point is replaced by a chord (arriving at
the support on either side like this: \psdiag{0.5}{3.5}{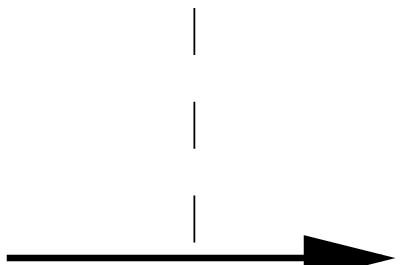}), we have adjusted the
sign in the exponent of the image of the
crossings.

With this modified version, nothing can keep us from forgetting the
orientations both of the 3-tangles and of the 3-diagrams. \hfill $\Box$

\section{The $(\frakg_2,V)$-weight system extended to 3-dia\-grams}

Let ${\cal A}$ be the subspace of the morphisms of ${\hat{\cal D}}_3$
generated by diagrams with support $S^1$ and without univalent
vertices.
It has been known for some time that given a Lie algebra $\frakg$, a
representation $V$ of $\frakg$, and an ad-invariant symmetric non-degenerate
bilinear form on $\frakg$, one can construct a linear
form---called weight system---on ${\cal A}$ if an orientation on the support
is added (see e.g.\ [BN 1], section 2.4,
or [V], section 6)\footnote{Our approach is closer to the one
described in [V].}. In this section, we construct an unoriented and
extended
version of the weight system one gets from \g with its
7-dimensional standard representation\footnote{``Standard'' in the
sense that every irreducible representation of \g occurs as a

direct summand of some tensor power of $V$.} and the bilinear
form $h\kappa$ (where $h$ is any non-zero complex number and
$\kappa$ the Killing form on \gf). We
do this by defining a functor $\Psi$ from \D \/ to ${\cal C}(\frakg_2)$,
a modification of a subcategory of the category of representations of
\gf.\\

Denote by $V$ the 7-dimensional irreducible representation and by $L$
the adjoint representation of \gf. Let the highest weights of
these representations be $(1,0)$ and $(0,1)$ respectively.

The following  fact assures that $V$ is selfdual (i.e.\ $V \cong V^*$) and
that there exist \gf-linear embeddings $i_\C$ and $i_V$ from $\C$
and $V$ respectively
into $V \otimes V$ unique up to scalars.
\begin{fact}\label{fact}
$$V \otimes V \cong {\bf C} \oplus V \oplus L \oplus W $$
$$\mbox{with}\quad
  \mbox{Sym}^2V  \cong {\bf C} \oplus W,\quad \wedge^2 V \cong V \oplus L$$
where  $W$ is the irreducible representation of \g with highest weight
$(2,0)$.
\end{fact}

For the construction of the functor $\Psi$, we will also need the following \gf-linear maps:
\begin{itemize}
\item $p_{\C} : V \otimes V\rightarrow \C$ and $p_V:V\otimes V \rightarrow
V$, the projections belonging to the embeddings $i_\C$ and $i_V$
(i.e.\ $p_\C \circ i_\C=id_\C,\, p_V\circ i_V=id_V$).
\item $flip_{X\otimes Y}: X\otimes Y \rightarrow Y\otimes X$, the \gf-linear
map taking $x\otimes y$ to $y\otimes x$  $(\forall x \in X, y\in Y$ for $X,Y \in
\{V,L\})$.
\item $cas: \C \rightarrow L\otimes L$ with
$(h\kappa)\circ cas = 14 id_\C$. Observe that $cas$ maps 1 to the
Casimir belonging to $h\kappa$.
\end{itemize}

Now we define the category ${\cal C}(\frakg_2)$ that will be the
target of the functor $\Psi$. Let ${\hat h}$ be a formal parameter.
\begin{df}
The category ${\cal C}(\frakg_2)$ is the monoidal
${\bf C}[\hspace{-0.1em}[\hat {h}]\hspace{-0.1em}]$-category with
objects {\em Obj(${\cal C}(\frakg_2)$)}$:=\{{\bf C}[\hspace{-0.1em}[\hat {h}]\hspace{-0.1em}]
\otimes_\C U\,|\, U \mbox{ is a tensor product over }\C \mbox{ with factors }
V \mbox{ and } L\}$ and with the following morphism spaces:\\ \centerline{$\mbox{{\em Mor}}_{{\cal C}(\frakig_2)}(
{\bf C}[\hspace{-0.13em}[\hat {h}]\hspace{-0.13em}]\otimes_\C U_1,
{\bf C}[\hspace{-0.13em}[\hat {h}]\hspace{-0.13em}]\otimes_\C U_2):=
{\bf C}[\hspace{-0.13em}[\hat
{h}]\hspace{-0.13em}]\otimes_\C\mbox{{\em
Hom}}_{\frakig_2}(U_1,U_2).$}
\end{df}

The definition of $\Psi$ is contained in the following proposition whose proof will be omitted, because
it just consists in checking straightforwardly that $\Psi$ respects all relations
required.
\begin{prop}
For any \gf-linear embeddings $i_\C$ and $i_V$ of $\C$ and $V$ respectively into
$V\otimes V$, there exists $k \in \C$ for which
we obtain a well defined $\C[\hspace{-0.1em}[\hat {h}]\hspace{-0.1em}]$-linear
monoidal functor $\Psi: \hat{\cal D}_3 \rightarrow {\cal C}(\frakg_2)$ by setting
\begin{itemize}
\item[(i)] $\Psi(\circ):= \C[\hspace{-0.1em}[\hat {h}]\hspace{-0.1em}]\otimes L\quad\quad\quad
\Psi(\bullet):= \C[\hspace{-0.1em}[\hat {h}]\hspace{-0.1em}]\otimes V$
\item[(ii)] $\Psi(\psdiag{2}{6}{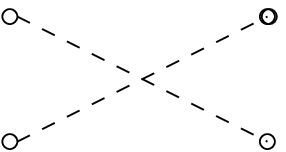}):= 1\otimes flip_{L\otimes L}\quad\hspace{3mm} \quad
\Psi(\psdiag{2}{6}{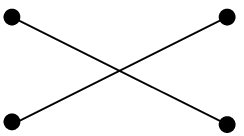}):= 1\otimes flip_{V\otimes V}$\newline
$\Psi(\psdiag{2}{6}{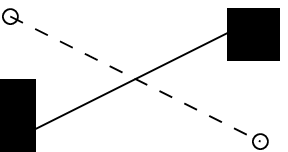}):= 1\otimes flip_{L\otimes V}\quad\quad\quad \Psi(
\psdiag{2}{6}{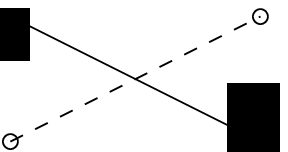}):= 1\otimes flip_{V\otimes L}$
\item[(iii)] $\Psi(\psdiag{2}{6}{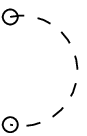}):=
1\otimes h\kappa,\quad\quad\quad\quad\quad\quad\quad\quad
\Psi(\psdiag{2}{6}{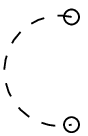}):= 1\otimes cas$\newline
$\Psi(\psdiag{2}{6}{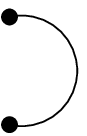}):=1\otimes 7p_\C\quad\quad\quad\quad\quad\quad\quad\quad \Psi(\psdiag{2}{6}{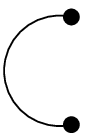}):= 1\otimes i_\C$
\item[(iv)] $\Psi(\psdiag{2}{6}{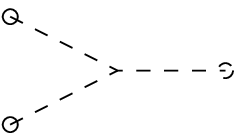}):={\hat h}\hspace{0.5ex}\otimes\hspace{0.5ex} $Lie bracket on $\frakg_2$\newline
$\Psi(\psdiag{2}{6}{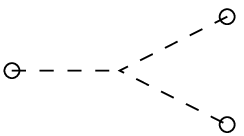}):={\hat h}\hspace{0.5ex}\otimes\hspace{0.5ex}$dual \footnote{Observe that $h\kappa$ induces an isomorphism from $L$ to $L^*$.}
of the Lie bracket on $\frakg_2$
\item[(v)] $\Psi(\psdiag{2}{6}{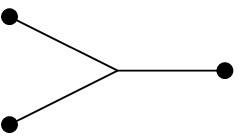}):= 1\otimes kp_V\quad\quad\quad\quad\quad\quad  \Psi(\psdiag{2}{6}{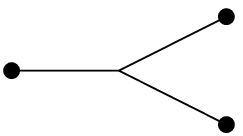})
:= 1\otimes i_V$
\item[(vi)]$\Psi(\psdiag{2}{6}{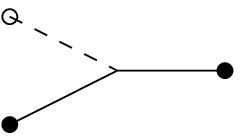}):={\hat h}\hspace{0.5ex}\otimes\hspace{0.5ex}$representation\newline
$\Psi(\psdiag{2}{6}{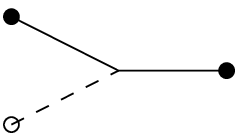}):=-{\hat h}\hspace{0.5ex}\otimes\hspace{0.5ex} $representation. \hfill$\Box$
\end{itemize}
\end{prop}

\begin{rmk} {\em The factor $k$ in $(v)$ depends on the embedding $i_V:V \rightarrow
V\otimes V$. The value of $k$ for a fixed embedding $i_V$ can be found by
solving the equation $(1\otimes (p_\C \otimes id_V))\circ(1\otimes (id_V\otimes i_V)) = 1\otimes kp_V$
representing the fact that $\Psi(\pstext{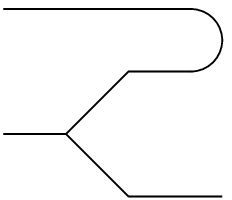}) = \Psi(\pstext{tril})$ must
hold.\newline
The factor 7 in $(iii)$ has been chosen to assure
$\Psi(\pstext{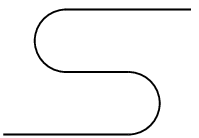})=\Psi(\pstext{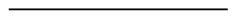})$; it is independent of the embedding
$i_\C$.\newline
The formal parameter $\hat h$ has been introduced to assure the
existence of $\Psi(D)$ for every morphism $D$ of $\hat {\cal
D}_3$: Due to $\hat h$, the morphism spaces of ${\cal C}(\frakg_2)$
can be regarded as graded spaces in the obvious way. This makes $\Psi$
a grade preserving functor that is well defined in every degree
and hence on the whole of $\hat {\cal D}_3$.}
\end{rmk}

\begin{rmk}
{\em In the introduction to this section, we mentioned the construction
of a weight system $\hat{\Psi}$ (out of Lie algebra information)
given in [BN 1] and [V]. In these papers, the support of the
diagrams is oriented, and so  the reader familiar with them might
ask if there is still a connection between $\Psi$ and $\hat{\Psi}$.
From the following observations, one can deduce that $\hat\Psi$ of an oriented
diagram of degree $m$ is exactly the degree $m$ part of $\Psi$ of the underlying
unoriented one:
\begin{tabbing}
\quad\quad\quad\=$1\otimes\hat{\Psi}( \pstext{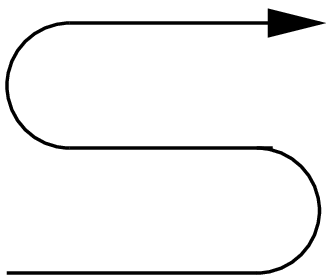})$\quad\= = \quad\=
${\Psi}( \pstext{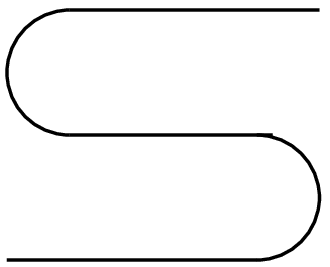})$ \quad\quad\quad\quad\=$\hat h^{\sharp\mbox{\scriptsize chords}}\otimes\hat{\Psi}(
\psdiag{2.5}{7.5}{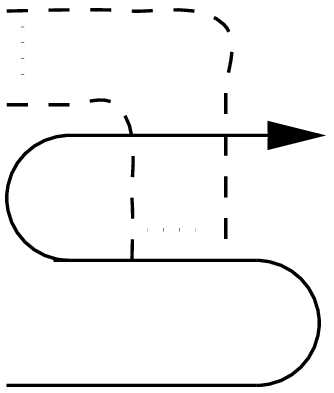})$\quad\= =\quad\=${\Psi}(
\psdiag{2.5}{7.5}{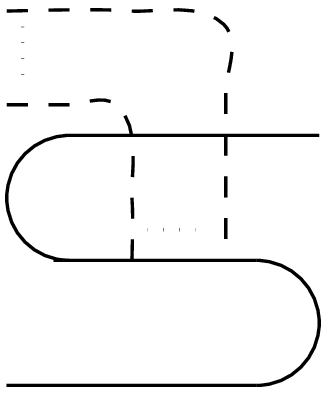})$\\[1em]
\>$1\otimes\hat{\Psi}( \pstext{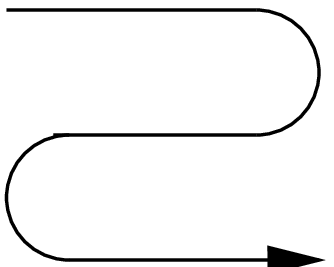})$\>=\>${\Psi}( \pstext{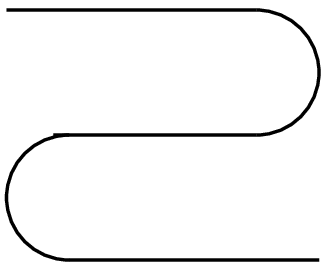})$
\>$\hat h^{\sharp\mbox{\scriptsize chords}}\otimes\hat{\Psi}( \psdiag{2.5}{7.5}{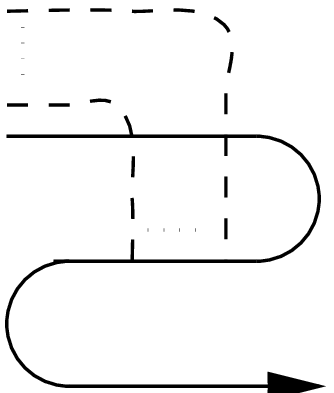})$\>=\>${\Psi}( \psdiag{2.5}{7.5}{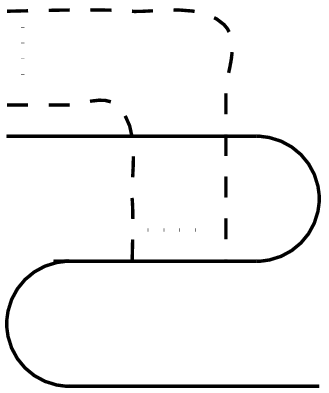})$
\end{tabbing}
(where ${\hat \Psi}(D)$ is obviously regarded as map
and not as tensor).}
\end{rmk}

\section{The skein relation for the $(\frakg_2,V)$-invariant}

What we have achieved so far, is the construction of an invariant for 3-tangles:
$\I:=\Psi\circ\hat Z_f$.
Unfortunately, we cannot evaluate it
directly because the expression known for the associator $\Phi$ is not explicit
enough to allow concrete computations.
But we will derive a skein relation, i.e.\ recursive rules by which
we can reduce the
problem to finding the values for planar 3-tangles (with these, we will
deal in the next section).
\par
The idea is to cut out a small neighbourhood of a crossing and insert
something else without changing the value of the invariant. The substitute
for the crossing has to be a linear combination of small, simple 3-tangles with
four univalent vertices. Obvious candidates for such are the inverse
crossing, \pstext{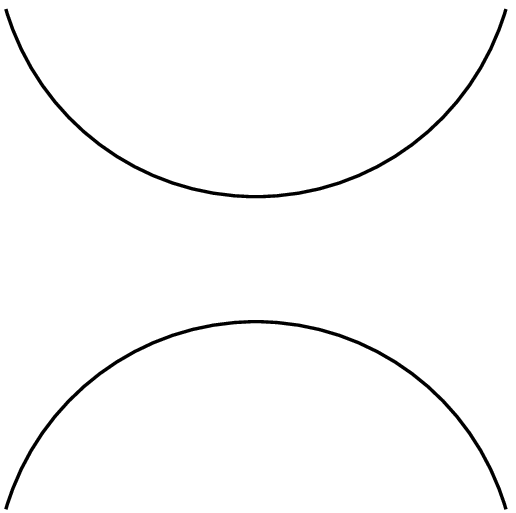}, and \pstext{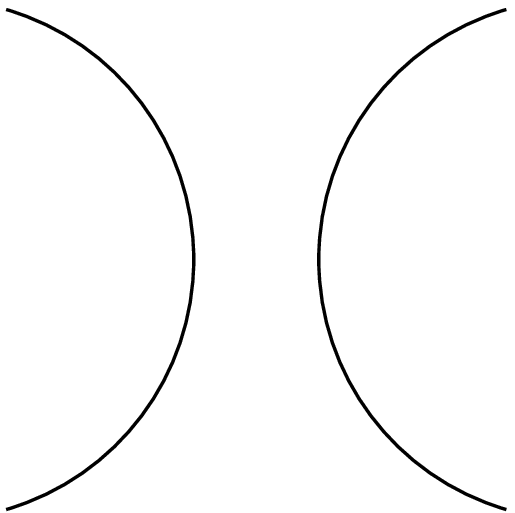};
as their values will prove to be linearly independent in the space of endomorphisms
of $V\otimes V$, these are not sufficient, and therefore, we
include \pstext{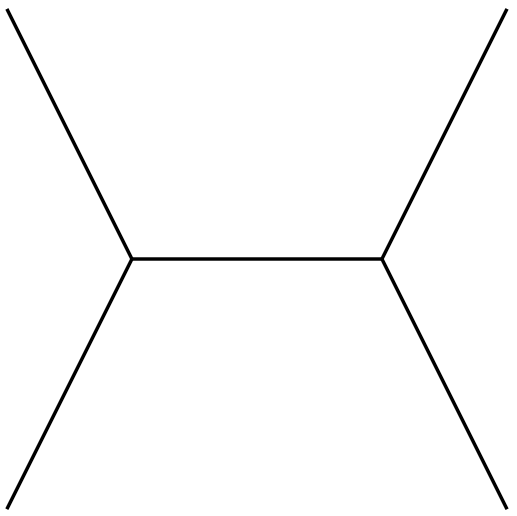} into our considerations.
\par
As $V\otimes V$ decomposes into four different irreducible representations
(namely {\bf C}, $V$, $L$ and $W$; see fact \ref{fact}), each $\frakg_2$-linear map
$V\otimes V \to V\otimes V$ is determined by the four corresponding
eigenvalues. To establish our skein relation, we have to ascertain the four eigenvalues of
\mbox{$\I$(\pstext{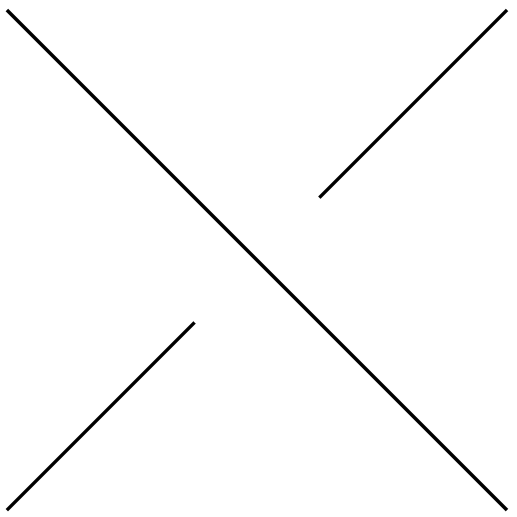}),} \mbox{$\I$(\pstext{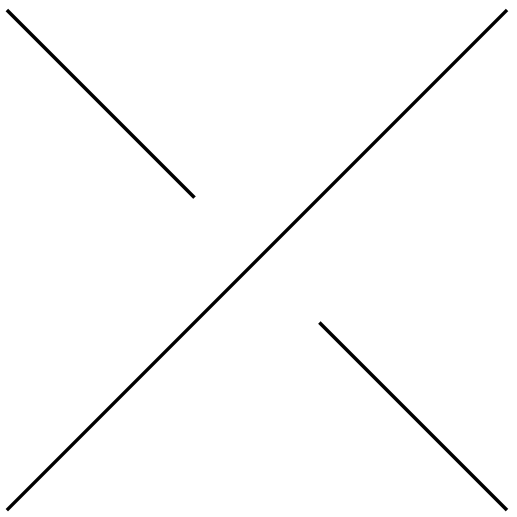}),}
\mbox{$\I$(\pstext{wbog}),}
\mbox{$\I$(\pstext{sbog}),} and $\I$(\pstext{h}).
\bigskip

{\bf Eigenvalues of $\I$(\pstext{pkreuz}) and
$\I$(\pstext{nkreuz})}
\\ The eigenvalues of $\I$(\pstext{pkreuz}) and
$\I$(\pstext{nkreuz})
are the products of the corresponding eigenvalues of
$\Psi(e^{\mp\frac{1}{2}\psdiag{0}{2}{hh}})$ and
$\Psi(\pstext{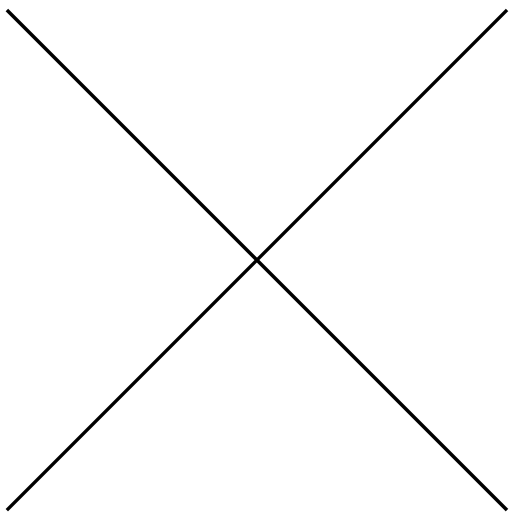})$.
\\ Quite a bit of explicit calculation yields that the eigenvalues of
$\Psi(\psdiag{0}{3}{hh})$ are
$\frac{\hat h^2}{2h}$ on {\bf C}, $\frac{\hat h^2}{4h}$ on $V$, 0 on $L$ and
$-\frac{\hat h^2}{12h}$ on $W$; accordingly, the eigenvalues of
$\Psi(e^{\mp\frac{1}{2}\psdiag{0}{2}{hh}})$
are $e^{\mp \frac{\hat h^2}{4h}}$, $e^{\mp \frac{\hat h^2}{8h}}$, 1, and
$e^{\pm \frac{\hat h^2}{24h}}$, respectively.
\\ The eigenvalue of $\Psi(\pstext{kreuz})$ is 1 on Sym$^2V={\bf C}\oplus W$ and -1
on $\wedge^2V=V\oplus L$.
\medskip

{\bf Eigenvalues of $\I$(\pstext{wbog})}
\\ As $\I$(\pstext{wbog}) is the identity on $V\otimes V$, its only
eigenvalue is 1.
\medskip

{\bf Eigenvalues of $\I$(\pstext{sbog})}
\\ Let $c$ be the scalar by which $\Psi(\pstext{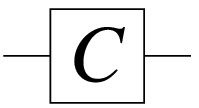})$ operates on $V$. Then we
have:
$$\I(\pstext{sbog}) = \Psi(\psdiag{2}{6}{Zim3}\,\psdiag{2}{6}{Zim7})
= c\Psi(\pstext{sbog})
= 7c\,i_{\bf C}\circ p_{\bf C}.$$
Hence $\I$(\pstext{sbog}) is 0 everywhere except on {\bf C}, where it has
the eigenvalue $7c$.

\begin{rmk}\label{unknot}
{\em Observe that:
$$\I(\pstext{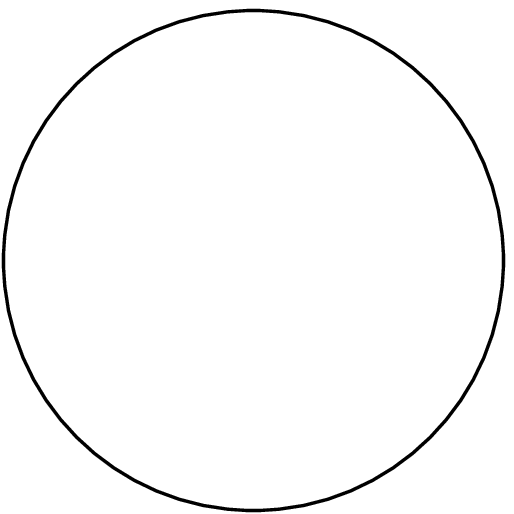}) = \Psi\circ\hat Z_f(\pstext{kreis})
=c\Psi(\pstext{kreis})=7c.$$
As Piunikhin has shown in [P] that for framed knots the Reshetikhin-Turaev
quantum invariants yield the same values as the invariants obtained by using
the corresponding weight systems, we can use a result of Rosso and
Jones in [RJ] to determine the value of $c$:
$$\I(\pstext{kreis})=\prod_{\alpha\in\Delta_+}\frac{e^{\frac{\hat h^2 h}{2}
(\lambda+\delta,\alpha)}-e^{-\frac{\hat h^2 h}{2}(\lambda+\delta,\alpha)}}{e^{\frac{\hat h^2 h}{2}
(\delta,\alpha)}-e^{-\frac{\hat h^2 h}{2}(\delta,\alpha)}} $$
\begin{tabbing}
where\quad\quad\=$\Delta_+$ is a possible choice for the set of positive roots
of $\frakg_2$\\
\>$\lambda$ is the highest weight of $V$\\
\>$\delta:=\sum_{R\in \Delta_+}R$\\
\>$(\,\,,\,)$ is the bilinear form on the weight space of $\frakg_2$ induced
by\\
\>the bilinear form $h\kappa$ on $\frakg_2$.
\end{tabbing}
Simplifying this expression and setting $q:=e^{\frac{\hat h^2}{24h}}$, we get:
$$7c=q^5+q^4+q+1+q^{-1}+q^{-4}+q^{-5}.$$
}
\end{rmk}

{\bf Eigenvalues of $\I$(\pstext{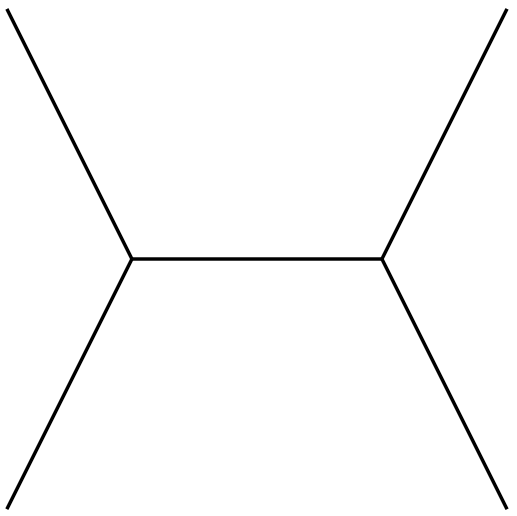})}
\\ Let $a$ (resp.\ $b$) be the eigenvalue of $\Psi(\pstext{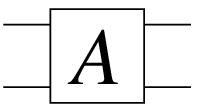})$ (resp.\
$\Psi(\pstext{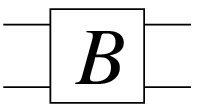})$) on $V$. Then we have:
$$\I(\pstext{H}) = \Psi(\psdiag{2}{6}{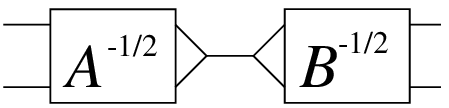}) =
\frac{\hat r^2}{\sqrt{ab}}\Psi(\pstext{H}).$$
Since $\Psi(\pstext{H})$ is 0 everywhere except on $V$, the parameter
$\hat r$ occurs nowhere but in the eigenvalue of $\I(\pstext{H})$
on $V$, which is $\frac{k\hat r^2}{\sqrt{ab}}$ (and $\ne 0$); therefore, we
do not have to care about the factor $\frac{1}{\sqrt{ab}}$, but can simply
shift the possibility of choice from $\hat r$ to
$r:=\frac{k\hat r^2}{\sqrt{ab}}$.

\begin{rmk}
{\em It is nonetheless possible to determine $ab$ by using the
following result of [LM 4] section 4:
$AB=(C^{-1}\otimes C^{-1})\Delta(C_1)$, and lemma \ref{lemma2}:
\begin{eqnarray*}
ab\Psi(\pstext{H}) & = & \Psi(\pstext{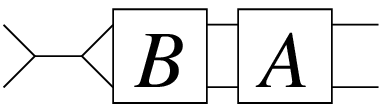})\\
                   & = & \Psi(\psdiag{3}{9}{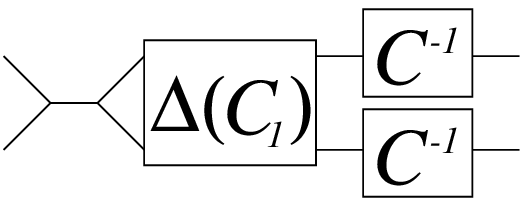})\\
                   & = & \Psi(\psdiag{3}{9}{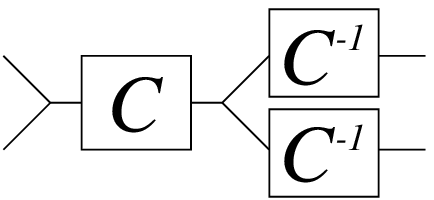})\\
                   & = & \frac{1}{c}\Psi(\pstext{H}),
\end{eqnarray*}
thus $$ ab = \frac{1}{c}.$$ }
\end{rmk}
\medskip

\begin{rmk}
{\em The invariant \It takes actually values in
${\bf C}[\hspace{-0.13em}[\hat h]\hspace{-0.13em}]$, but as long as
we do not want to fix $h$ and $r$, we can regard them as elements
of ${\bf C}[\hspace{-0.13em}[\hat h,\frac{1}{h},
r]\hspace{-0.13em}]$ and, accordingly, ${\cal C}(\frakg_2)$ as
${\bf C}[\hspace{-0.13em}[\hat h,\frac{1}{h},r]\hspace{-0.13em}]$-category.}
\end{rmk}

To summarize (recall that $q=e^{\frac{\hat h^2}{24h}}$):
\vspace{5mm} \\
\begin{tabular}{|c||c|c|c|c|c|} \hline
$\begin{array}{c} \mbox{{\footnotesize Eigen-}} \\
\mbox{{\footnotesize value on}}\end{array}$
& $\I$(\pstext{pkreuz})
& $\I$(\pstext{nkreuz})
& $\I$(\pstext{wbog}) & $\I$(\pstext{sbog})
& $\I$(\pstext{H}) \\ \hline\hline \rule{0cm}{3ex}
{\bf C} & $q^{-6}$ & $q^6$
& 1 & $7c$ & 0 \\
$V$ & $-q^{-3}$ & $-q^3$
& 1 & 0 & $r$ \\
$L$ & $-1$ & $-1$ & 1 & 0 & 0 \\
$W$ & $q$ & $q^{-1}$ & 1 & 0 & 0 \\ \hline
\end{tabular}
\vspace{5mm}

The leftmost column can be expressed as a linear combination of the other
columns; i.e.\ substituting this linear combination of the 3-tangles
\pstext{nkreuz},
\pstext{wbog}, \pstext{sbog}, and \pstext{H}
 for a crossing \pstext{pkreuz} in a 3-tangle does not
change the value of $\I$.

\begin{thm}\label{thm4}
For the invariant $\I$, the following skein relation holds:
$$\I(\pstext{pkreuz}) = \alpha \I(\pstext{nkreuz})+
\beta \I(\pstext{wbog})
+\gamma \I(\pstext{sbog}) + \delta \I(\pstext{H})$$

\begin{tabbing}
where xx\=$\lambda$ x \= := x \= $\frac{\alpha\gamma+\beta}{1-\alpha^2}$\kill
where \>$\alpha$ \> := \> $q $ \\
      \>$\beta$  \> := \> $q-1$ \\
      \>$\gamma$ \> := \> $\frac{1}{7c}(-q^7+q^{-6}-q+1)$ \\
      \>$\delta$ \> := \> $\frac{1}{r}(q^4-q^{-3}-q+1)$.
\end{tabbing}
\schluss
\end{thm}

Since $\I$ is a monoidal functor and invariant under ambient isotopy,
we can deduce another skein relation as follows:
\begin{eqnarray*}
\I(\pstext{nkreuz}) & = & \I(\pstext{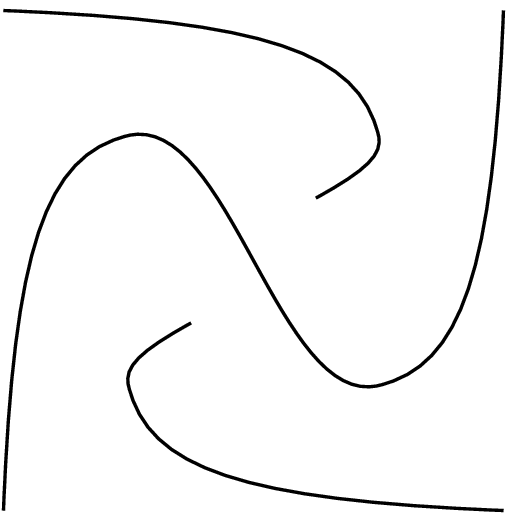}) \\
& = &  \alpha \I(\pstext{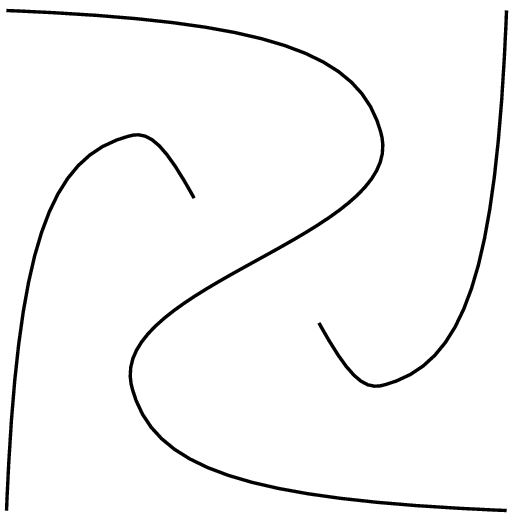})+ \beta \I(\pstext{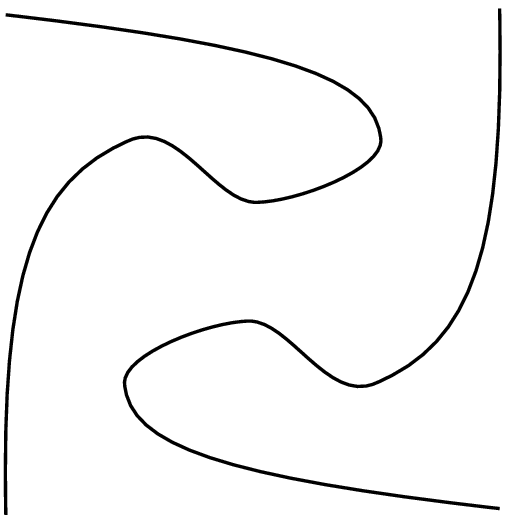})
+\gamma \I(\pstext{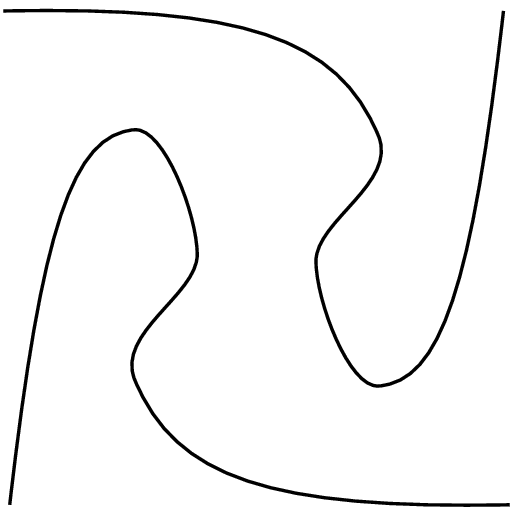}) + \delta \I(\pstext{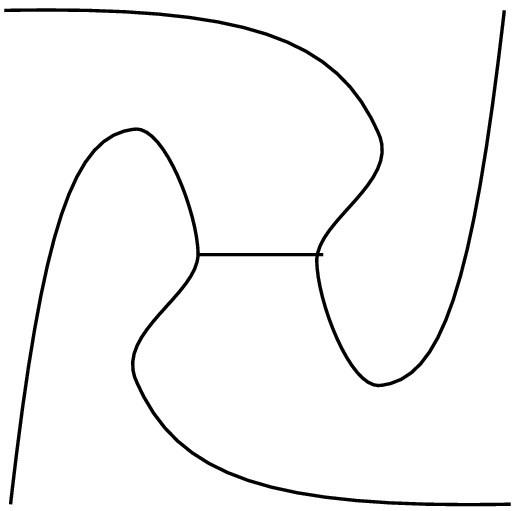}) \\
& = &  \alpha \I(\pstext{pkreuz})+ \beta \I(\pstext{sbog})
+\gamma \I(\pstext{wbog}) + \delta \I(\pstext{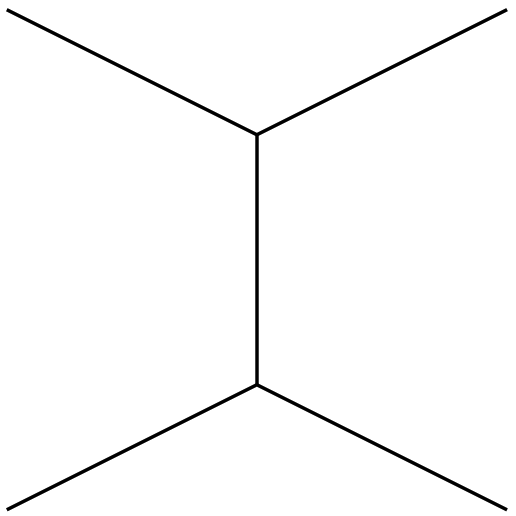}).
\end{eqnarray*}

Combining these two versions of the skein relation, we obtain:

\begin{cor}
For the invariant $\I$, the following skein relation holds
\footnote{Setting
$r=-(q^2+q+1+q^{-2}+q^{-3}+q^{-4})$,
we obtain the skein relation given in [K].}:
$$\I(\pstext{pkreuz}) = \lambda \I(\pstext{wbog})+
\mu \I(\pstext{sbog})
+\rho \I(\pstext{H}) + \sigma \I(\pstext{I})$$

\begin{tabbing}
where xx\=$\lambda$ x \= := x \= $\frac{\alpha\gamma+\beta}{1-\alpha^2}$\kill
where   \>$\lambda$ \> :=  \> $\frac{\alpha\gamma+\beta}{1-\alpha^2}$ \\
        \>$\mu$ \> := \> $\frac{\alpha\beta+\gamma}{1-\alpha^2}$ \\
        \>$\rho$ \> := \> $\frac{\delta}{1-\alpha^2}$ \\
        \>$\sigma$ \> := \> $\frac{\alpha\delta}{1-\alpha^2}.$
\end{tabbing}
\schluss
\end{cor}

It is clear that by means of this relation, every 3-tangle can be reduced
to a linear combination of planar 3-tangles.

\begin{rmk}
{\em The invariant $\I = \Psi\circ\hat Z_f$ of closed oriented 3-nets itself is
not a Vassiliev invariant, but ``consists of'' Vassiliev invariants in the
following sense: For each $m \in ${\bf N}, let $\I^{(m)}$ be the
function with values in ${\bf C}[\hspace{-0.13em}[\frac{1}{h},
r]\hspace{-0.13em}]$ such that
$\I = \sum_{m=0}^{\infty} \I^{(m)}\hat h^{2m}$. Then $\I^{(m)}$ is
an invariant of type $m$.}
\end{rmk}


\section{Values of the $(\frakg_2,V)$-invariant on closed planar 3-nets}
In this section, we show how the value of \It on a closed planar 3-net can
be calculated recursively.

The following lemma assures that it is sufficient to consider connected
3-nets.

\begin{lem}
If a closed 3-net $N$ is equivalent to a 3-net consisting of two closed
3-nets $N_1$ and $N_2$ with $N_1 \subset \R^-\times \R^2, N_2\subset
\R^+\times \R^2$, then
$$\I(N) = \I(N_1)\cdot \I(N_2).$$    \hfill $\Box$
\end{lem}

As, by definition, every planar 3-net is  equivalent to a 3-net contained in
$\R^2\times \{0\}$ with only upward pointing vectors assigned, it is enough
to calculate \It for these. Therefore, we assume that all planar 3-nets in this
section are of this type.

\begin{df}
Let $N$ be a planar 3-net. A {\em mesh} of $N$ is the closure of a
connected component of $(\R^2\times \{0\})\backslash N$. A {\em $n$-mesh}
is a mesh with $n$ trivalent vertices in the boundary.
\end{df}

We will show how in any non-empty connected closed planar 3-net the number
of meshes can be reduced without changing the value of the invariant. As
we know that $\I$(empty 3-net) = 1 (the empty 3-net is the unity in \topf
\hspace{0mm} and $\I$ is a monoidal functor),
this will allow us to calculate the
invariant of a closed planar 3-net recursively.

\begin{prop}
Let $N$ be a non-empty connected planar 3-net with $m$ meshes. Then there
exist closed planar 3-nets $N_1,\ldots , N_r$  (not necessarily connected)
with fewer than $m$ meshes and coefficients $\lambda_1,\ldots,\lambda_r \in
\C$ such that $\I(N)=\sum_{i=1}^r\lambda_i\I(N_i)$.
\end{prop}

{\bf Proof:} The idea is to cut out a mesh and replace it by a linear
combination of pieces that lead to 3-nets with fewer meshes.

Thanks to the following lemma, the mesh we want to cut out can always be
chosen to be a simply connected $n$-mesh with $n\leq 5$.

\begin{lem}\label{km}
Let $N$ be a planar non-empty closed connected 3-net. Then $N$ has at least
one simply connected $n$-mesh with $n\leq 5$.
\end{lem}

{\bf Proof of the lemma:} If $N$ is an embedded $S^1$, $N$ has a bounded
0-mesh, and so in this case, the lemma holds.\newline
Let $\sharp$ denote ``number of\," and let $N$ be a 3-net without 0-mesh.

\underline{Observation 1}: $\sharp$ vertices of $N - \sharp$ edges of $N +
\sharp$ meshes of $N =2$ (Euler characteristic of $S^2$).

\underline{Observation 2}: If we assign to each mesh $M$ the appropriate
part of the contributions of its vertices, its edges, and  its region, i.e.\
$$\chi_M:=\left\{\begin{array}{l}\quad\frac{1}{3}\: \sharp\mbox{ vertices of
}M\\ +\frac{1}{3}\: \sharp\mbox{ vertices of }M\mbox{ for  which all adjacent
edges belong to } M\\ +\frac{1}{3}\: \sharp\mbox{ vertices of }M\mbox{ that
belong only to }M\\ -\frac{1}{2}\: \sharp\mbox{ edges of }M\\
 -\frac{1}{2}\: \sharp\mbox{ edges of }M\mbox{ that belong only to } M\\
 +1,
\end{array}\right.$$ $$\mbox{then } \sum_{M \mbox{ {\scriptsize mesh of} } N}\chi_M =
2.$$

\underline{Observation 3}: For a simply connected $n$-mesh $M$ of $N$, we
have $\chi_M= 1-\frac{1}{6}n$.

\underline{Observation 4}: If the unbounded mesh $M'$ of $N$ does not
contain an edge belonging to $M'$ only, then $\chi_{M'}\leq 1$.

If $N$ does not contain an edge that belongs to one mesh only, the lemma
is a consequence of observations 2, 3 and 4\footnote{Note that any bounded mesh $M$ that does not
contain an edge belonging only to $M$ is simply connected.}.

Now suppose that there are edges $e_1,\ldots, e_k$ such that $e_j$ belongs
 to only one mesh $M_j$. Note that forgetting such an edge $e_j$ would split
 the 3-net $N$ into two connected components $N_{j1}$ and $N_{j2}$. As there are only
finitely many edges $e_j$, there is  an edge $e_i$ for which at least for one
$k\in \{1,2\}$ $N_{ik}$ satisfies the following properties:\newline
($i$)$\,\quad$ $N_{ik}$ does not contain an edge that belongs to only one mesh;\newline
($ii$)\quad any mesh $M$ of $N$ that is also a mesh of $N_{ik}$ is bounded.\newline
The sum $\sum \chi_M$ over all meshes $M$ of $N$ mentioned in ($ii$) is greater than
0 (look at $\sum \chi_M$ with $M$ considered as mesh of $N_{ik}$, use observation 4, and subtract
$\frac{1}{6}$ for the influence of $e_i$),
and thus observation 3 guarantees that at least one of
these meshes has fewer than 6 vertices.\hfill $\Box$ of lemma.\\

To replace $n$-meshes for $n\leq 5$, we will use the following lemma:

\begin{lem}
The following equations hold:
\begin{tabbing}
(o)\quad \=$\I (\psdiag{1}{3}{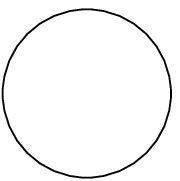})$\quad\quad\= =\quad\quad \=$7c\I(\quad\quad)=7c$,\\[5pt]
(i)\>$\I(\psdiag{1}{3}{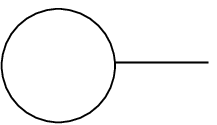})$\>=\> 0,\\[5pt]
(ii)\>$\I(\psdiag{1}{3}{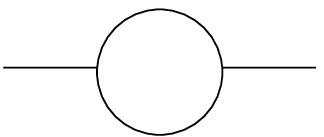})$\>=\>$r\I(\pstext{ek33})$,\\[5pt]
(iii)\>$\I(\pstext{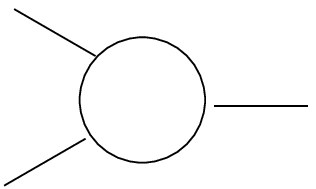})$\>=\>$t\cdot\I(\pstext{tril})$ with
$t:=\frac{1}{\delta}(-q^3+\alpha q^{-3}-\gamma)$,\\[5pt]
(iv)\>$\I(\pstext{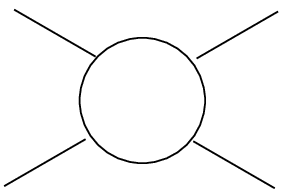})$\>=\>$\frac{r^2q^5}{g(q^4+1)}(\I(\pstext{wbog})
+\I(\pstext{sbog}))$\\[5pt]
\>\>\>$+\frac{rq^2(q^2+1)}{g}(\I(\pstext{h})+
\I(\pstext{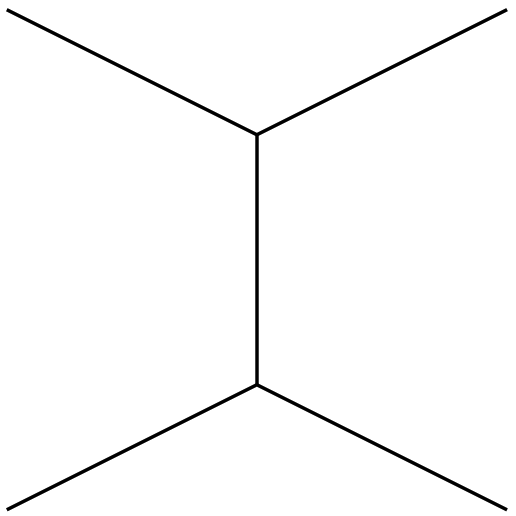}))$\\[0em]
\>\>\>with $g:=q^6+q^5+q^4+q^2+q+1$, \\[5pt]
(v)\>$\I(\psdiag{3}{9}{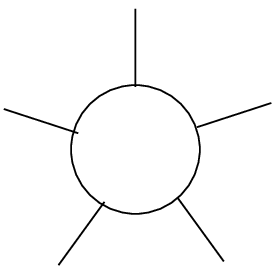})$\>=\>$-d(\I({\mbox{\psdiag{3}{9}{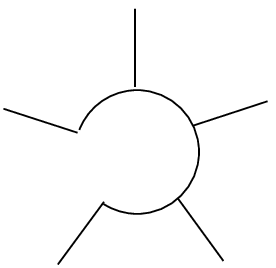}}})+
\I({\mbox{\psdiag{3}{9}{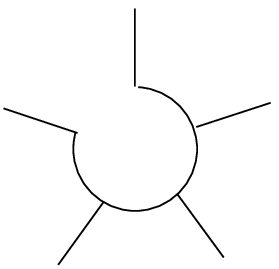}}})+\I({\mbox{\psdiag{3}{9}{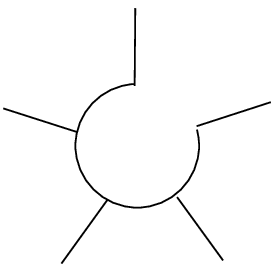}}})$\\[5pt]
\>\>\>$\quad\quad+\I({\mbox{\psdiag{3}{9}{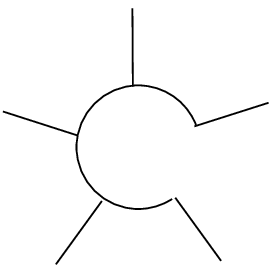}}})+
\I({\mbox{\psdiag{3}{9}{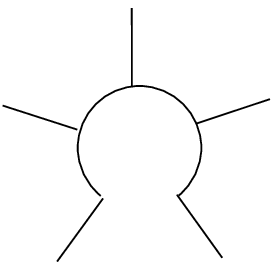}}}))$ \\[5pt]
\>\>\>$-d^2(\I({\mbox{\psdiag{3}{9}{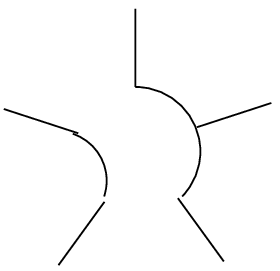}}})+ \I({\mbox{\psdiag{3}{9}{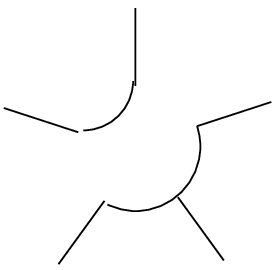}}})+
\I({\mbox{\psdiag{3}{9}{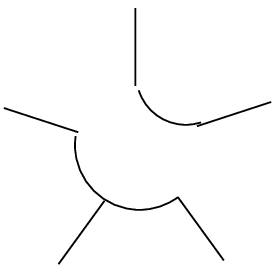}}})$\\[5pt]
\>\>\>$\quad\quad+\I({\mbox{\psdiag{3}{9}{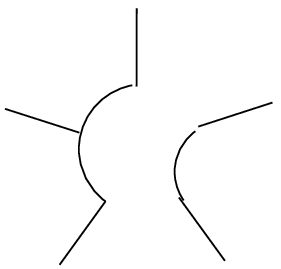}}})+
\I({\mbox{\psdiag{3}{9}{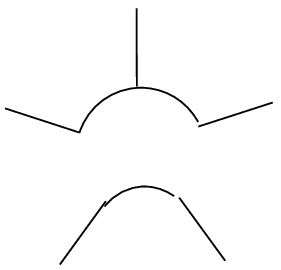}}}))$ \\[5pt]
\>\>\>with $d:=\frac{rq^3}{g}$.
\end{tabbing}
\end{lem}

{\bf Proof of the Lemma:} Equation ($o$) is proved in remark
\ref{unknot}.

Equation ($i$) is true because $\I(\psdiag{1}{3}{1mesh})$ is a $\frakg_2$-linear
map from $\C$ to $V$ and must therefore be 0.\newline
Equation ($ii$), we get from
$\I(\pstext{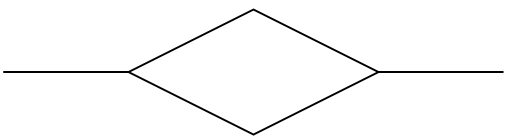})=\Psi(\pstext{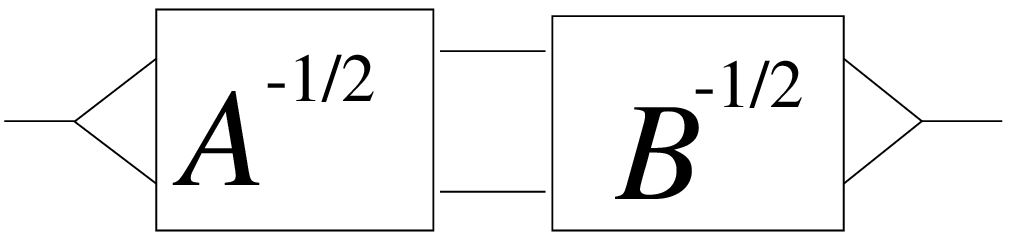})\stackrel{\mbox{{\scriptsize Sec.\ 4}}}{=}
r\Psi(\pstext{ek33})$.\newline
To get equation ($iii$), we use the skein relation given in theorem \ref{thm4}
(rotated by $90^\circ$):
$$\I(\psdiag{1}{3}{3mesh})=\frac{1}{\delta}(\I(\pst{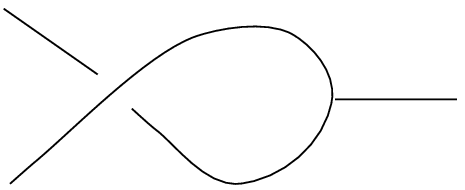})-\alpha\I(
\pst{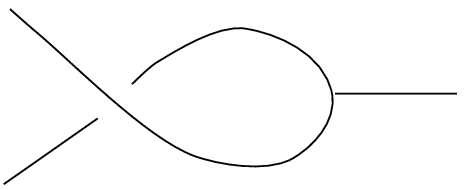})-\beta\I(\pst{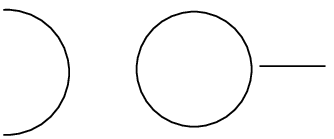})-\gamma\I(\pst{tril}).$$
Note that
$\I(\pst{5k13})=\I(\pst{tril})\circ\I(\pst{nkreuz})=-q^3\I(\pst{tril})$
(only the eigenvalue of $\I(\pst{nkreuz})$ on $V$ matters).\newline
Equation ($iv$) and ($v$) can be obtained in a similar way. For ($v$), it may
help to use that $\pstext{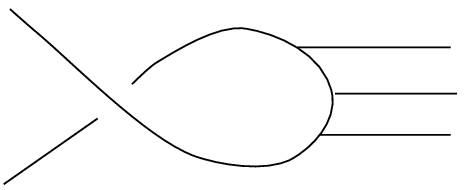}\sim\pstext{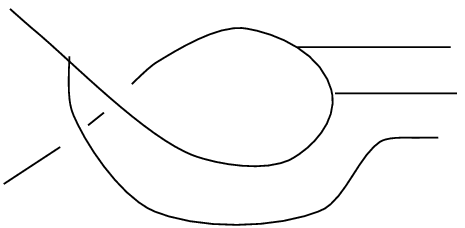}$.\hfill
$\Box$of lemma.\\

As \It is a monoidal functor, equations ($o$)-($v$) will still hold if the
3-nets depicted in the arguments of \It are parts of  bigger 3-nets that
are identical outside the depicted region for all arguments in the same
equation. The observation that for all equations, the 3-nets appearing on
the right-hand side have fewer meshes than the one on the left-hand side
concludes the proof of the proposition.\hspace*{\fill} $\Box$\\

Now we are able to compute the value of the $\I$-invariant for every
oriented closed planar 3-net recursively. Planar 3-tangles that are not closed can
often be reduced by the same technique, but because in this case, lemma
  \ref{km} is no longer true, it may happen that we get stuck before reaching
the empty 3-net (example:\pstext{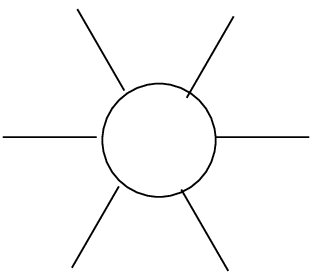}).\\

\section{Some examples}

To do explicit calculations, the following lemmas may be helpful, the first
comparing 3-nets to their mirror
images, the second suggesting some short cuts.

\begin{lem}
As taking mirror images essentially comes to changing crossings, the value of
$\I$ on the mirror image of a 3-net N is obtained by substituting
$q^{-1}$ for q in $\I(N)$.  \hfill $\Box$
\end{lem}

\begin{lem} \ \vspace{-5mm} \\
\begin{tabbing}
$\I(\psdiag{1}{4}{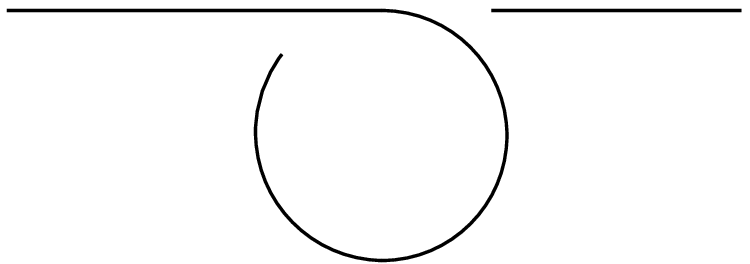})$ \= = \=
$q^6 \I(\psdiag{2.5}{6}{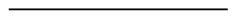})$\kill
$\I(\psdiag{1}{4}{fram1})$ \> = \>
$q^6 \I(\psdiag{2}{6}{strich})$\\[1ex]
$\I(\pstext{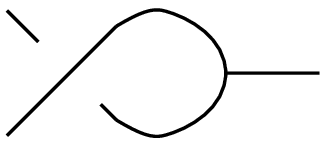})$ \> = \>
$-q^3 \I(\pstext{Zarg4})$\\[1ex]
$\I(\psdiag{1}{4}{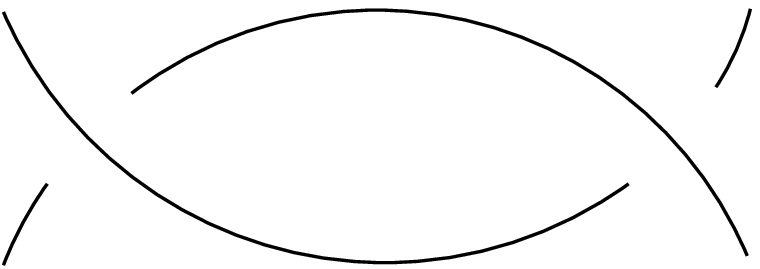})$ \> = \> $\beta \I(\pstext{pkreuz})$
$+\alpha \I(\pstext{wbog})$
$+\gamma q^{-6}\I(\pstext{sbog})$
$-\delta q^{-3} \I(\pstext{H})$
\end{tabbing}
\hfill $\Box$
\end{lem}

\begin{ex} \ \vspace{2mm} \\
\begin{tabular}{@{} l @{ = } p{11.5cm}}
$\I(\psdiag{1.5}{6}{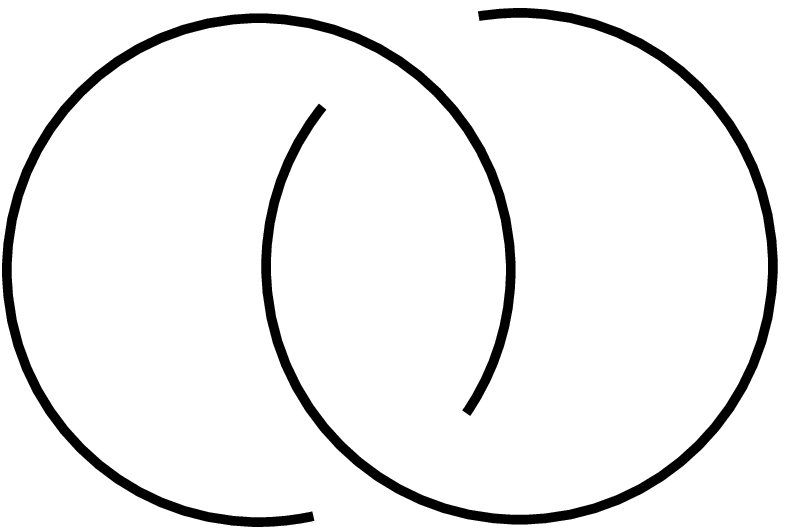})$ & $7c(q^7+q^5+q^2+1+q^{-2}+q^{-5}+q^{-7})$\\
[1ex]
$\I(\psdiag{1.5}{6}{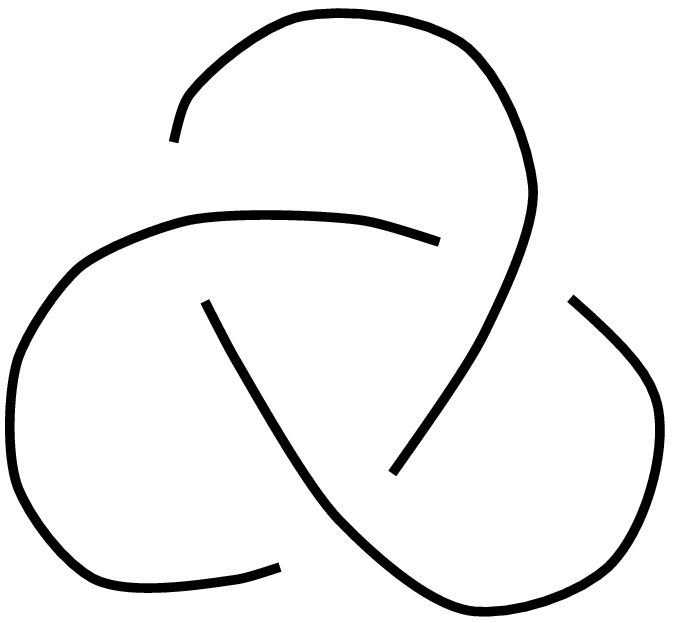})$ &
$7c(q^8+q^6-q^5+q^3-q^2+q-1+q^{-1}+q^{-4}-$
$2q^{-5}+2q^{-6}-q^{-7}-q^{-9}-q^{-10}+q^{-11}-q^{-12}+q^{-13})$\\ [1ex]
$\I(\psdiag{1.5}{6}{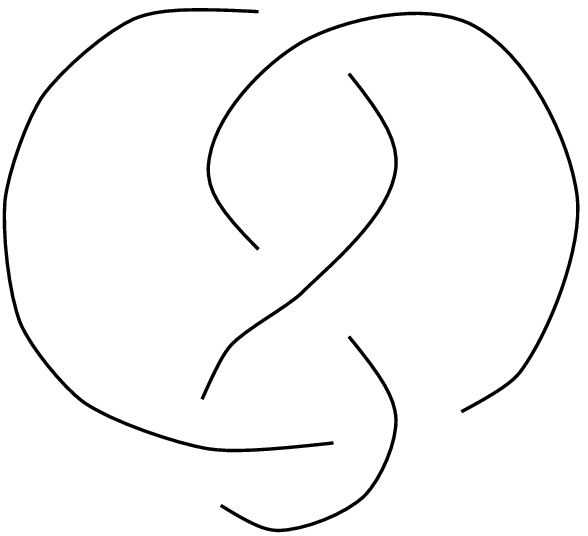})$ &
$7c(q^{14}-q^{13}+2q^{12}-2q^{11}+q^9-2q^8+4q^7-4q^6+4q^5-2q^4-q^3+3q^2-5q+5$
$-5q^{-1}+3q^{-2}-q^{-3}-2q^{-4}+4q^{-5}-4q^{-6}+4q^{-7}-2q^{-8}+q^{-9}-2q^{-11}$
$+2q^{-12}-q^{-13}+q^{-14})$\\ [1ex]
$\I(\psdiag{1.5}{6}{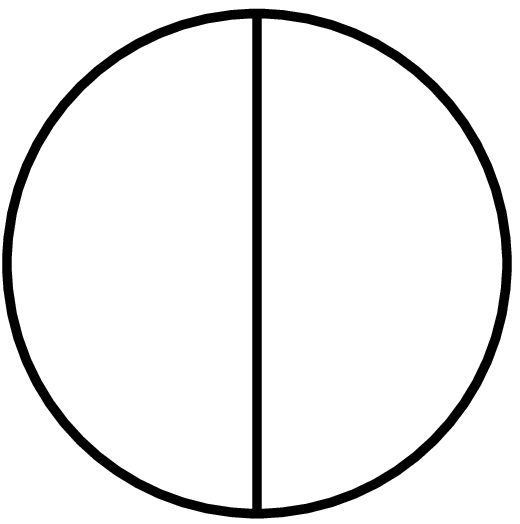})$ & $7cr$\\ [1ex]
$\I(\psdiag{1.5}{6}{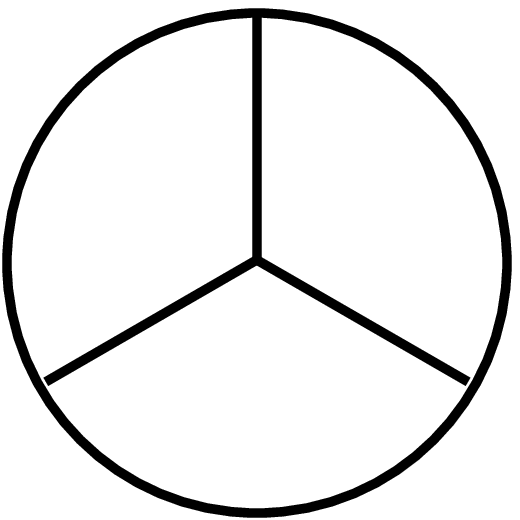})$ &
$-7cq\,\frac{q^2-q+1}{q^4+1}r^2$
\end{tabular}
\end{ex}



\newpage
\section{References}

\begin{tabular}{@{} l p{13cm}}
[BN 1] & D. Bar-Natan, {\em On the Vassiliev knot invariants}, Topology, 34
(1995), 423-472. \\
{[}BN 2] & D. Bar-Natan, {\em Non-associative tangles}, in {\em Geometric
topology}
(proceedings of the Georgia International Topology Conference), (W.H. Kazez, ed.),
139-183, Amer.\ Math.\ Soc.\ and International Press, Providence, 1997.\\
{[}FH]  & W. Fulton and J. Harris, {\sl Representation theory}, Graduate Texts in
Mathematics \# 129, Springer-Verlag 1991.\\
{[}H]   & J.E. Humphreys, {\sl Introduction to Lie algebras and representation
theory}, Graduate Texts in Mathematics \# 9, Springer Verlag 1994.\\
{[}K]   & G. Kuperberg, {\em The quantum $G_2$ link invariant}, Int.\ J. Math., 5
(1994), 61-85.\\
{[}LM 1] & T.T.Q. Le and J. Murakami, {\em Kontsevich's integral for the
Homfly polynomial and relations between values of multiple zeta functions},
Topology Appl., 62 (1995), 193-206.\\
{[}LM 2] & T.T.Q. Le and J. Murakami, {\em Kontsevich's integral for the
Kauffman polynomial}, Nagoya Math.\ J., 142 (1996), 39-65.\\
{[}LM 3]& T.T.Q. Le and J. Murakami, {\em The universal Vassiliev-Kontsevich
invariant for framed oriented links}, Compositio Math., 102 (1996), 41-64.\\
{[}LM 4] & T.T.Q. Le and J. Murakami, {\em Parallel version of the universal
Vassiliev-Kontsevich invariant}, J. pure appl.\ algebra, to appear.\\
{[}LMMO] & T.T.Q. Le, H. Murakami, J. Murakami and T. Ohtsuki: {\em A
three-manifold invariant via the Kontsevich integral}, Osaka
J. Math., to appear.\\
{[}MO]  & J. Murakami and T. Ohtsuki, {\em Topological quantum field theory for
the universal quantum invariant}, Commun.\ Math.\ Phys., to appear.\\
{[}P]   & S. Piunikhin, {\em Weights of Feynman diagrams, link polynomials, and
Vassiliev invariants}, J. knot theory ramifications, 4 (1995), 163-188.\\
{[}RJ]  & M. Rosso and V. Jones, {\em On the invariants of torus knots derived
from quantum groups},  2 (1993), 97-112.\\
{[}St]   & T. Stanford, {\em The functoriality of Vassiliev-type invariants
of links, braids, and knotted graphs}, J. knot theory ramifications, 3(3)
(1994), 247-262.\\
{[}V]   & P. Vogel, {\em Algebraic structures on modules of diagrams}, Invent.\ Math.,
to appear.
\end{tabular}

\end{document}